\documentclass[11pt]{article}
\usepackage{amsmath, amssymb,ytableau,young,bm,stackrel,textcomp,amsxtra, graphicx}
\usepackage{amsmath, amssymb,amsxtra, graphicx, float,tikz,ytableau,young,bm,stackrel,textcomp,}
\usepackage{bold-extra}
\usepackage{courier}
\usepackage{color}
\usepackage{tabularx}
\usepackage{tikz,pgf}
\newtheorem{thm}{{\sc Theorem}}[section]
\newtheorem{prop}[thm]{{\sc Proposition}}
\newtheorem{cor}[thm]{{\sc Corollary}}
\newtheorem{lem}[thm]{{\sc Lemma}}
\newtheorem{df}[thm]{{\sc Definition}}
\newenvironment{proof}{\begin{sc}\noindent Proof: \end{sc}}{
     \hbox to 2em{}\nobreak\hfill$\blacksquare$\par\medskip}
\newcommand{\LR}{\hbox{Little\-wood-Richard\-son}}
\newcommand {\PTab}{\hbox{PTab}}
\newcommand {\Pf}{\text{Perf\,}}
\newcommand {\SSYT}{\text{SSYT}}
\newcommand{\comm}[1]{}
\newcommand {\mmin}{\text{min}}
\newcommand {\mmax}{\text{max}}
\newcommand {\Rot}{\text{Rot}}
\newcommand {\NULL}{\text{NULL}}
\newcommand {\Evac}{\text{Evac}}
\newcommand{\Lus}{\text{Lus}}
\newcommand {\inv}{\text{inv}}
\newcommand{\wt}{\text{wt}}

\begin{document}

oooo 
\begin{center} {\bf Perforated Tableaux: A Combinatorial Model for Crystal Graphs in Type $A_{n-1}$}\end{center}

\medskip

\noindent
Glenn D. Appleby, Tamsen Whitehead*\\
{\em Department of Mathematics\\
and Computer Science,\\
Santa Clara University\\
Santa Clara,  CA USA}\\
*Corresponding author email: tmcginley@scu.edu,\\
gappleby@scu.edu\\

\noindent \mbox{} \hrulefill \mbox{}
\begin{abstract}
We present a combinatorial model, called \emph{perforated tableaux}, to study $A_{n-1}$ crystals, unifying several previously studied
combinatorial models. We identify nodes in the $k$-fold tensor product of the standard crystal with length $k$ words in $[n]= \{ 1, \ldots n\}$.
 We model this crystal with perforated tableaux (ptableaux), extending this identification isomorphically to biwords, RSK $(P,Q)$ tableaux pairs, and matrix models. In the ptableaux setting, crystal operators are more simply defined and we can identify highest weights visually without computation. We generalize the tensor products in the Littlewood-Richardson rule to all of $[n]^{\otimes k}$, and not just the irreducible crystals whose reading words come from semistandard Young tableaux.  We relate evacuation (Lusztig involution) to products of ptableaux crystal operators, and find a combinatorial algorithm to compute commutators of highest weight ptableaux.
\end{abstract}
\noindent \mbox{} \hrulefill \mbox{}

\noindent
{\bf Declarations:}\\
Funding: Not applicable\\
Conflicts of interest/competing interests: Not applicable\\
Availability of data and material: Not applicable\\
Code availability: Not applicable\\

\section{Introduction}
We introduce a new combinatorial model for the $GL_n$ crystal structure (of type $A_{n-1}$), called \emph{perforated tableaux} (or simply, \emph{ptableaux}), and study its relation to crystal structures on $[n]^{\otimes k}$, the $k$-fold tensor product of the standard $GL_n$ crystal, denoted by $[n]$, along with its relation to other combintorial models. Perforated tableaux are rectangular arrays containing concatenated horizontal strips of positive integers, and blank boxes (``perforations"), such that the content in non-blank entries weakly increases from left to right and strictly increases from top to bottom.  Entire rows of blanks are allowed; entire columns of blanks are not. For example the diagram below is a perforated tableau:

\ytableausetup{smalltableaux}
 \ytableausetup{centertableaux}
\[ \begin{ytableau}
\ & & 1 &2&5&6&6&6&7&8&8 \\
1&1&4&5&7&&8&8&&&\\
3&4&&7&8& 8&9&9&9&&\\
4&5&5&&&&&&&& \end{ytableau}\, . \]

While the $GL_n$ crystal structure in type $A_{n-1}$ is certainly well-studied, the perforated tableau model allows us to unify several previously separate models. The term ``perforated tableau" appears in the work of Benkart, Sottile, and Stroomer (see~\cite{Sottile}) in their work on tableau-switching algorithms. Our definition is similar to theirs, but not identical, and our use of perforated tableaux as models for crystal graphs is new.

\medskip

\medskip

The authors would like to thank Florence Maas-Gariepy for her very careful reading of an early version of this paper and her many useful suggestions, along with the contributions of the referee for their valuable suggestions and for bringing to our attention several results related to ours.

\subsection{Combinatorial Models for Crystals}

A \emph{combinatorial model} for a crystal graph is a set of objects used as nodes of the graph (here, of type $A_{n-1}$) so that we can define crystal operators $e_{i},f_j$ on such objects by combinatorial means. Often, one can then compute various multiplicities by enumerating some set of combinatorial objects; the overall goal of such models is to gain some insight into the operations involved, and potentially to simplify computations.

Below, we show how to obtain easily-defined bijections between perforated tableau and other well-studied combinatorial models, and describe why the ptableau model may have advantages over them.

\medskip

Bump and Schilling~\cite{BumpSchilling} define ${\cal B}_{n}$ as the \emph{standard} $GL_n$ crystal of type $A_{n-1}$ whose graph is composed of tableaux of one box:
\[ \framebox{1} \xrightarrow{f_{1}} \framebox{2} \xrightarrow{f_{2}} \cdots \xrightarrow{f_{n-1}} \framebox{n}. \]

Instead of ${\cal B}_{n}$, we'll let $[n]=\{1,2,\ldots , n \}$ denote both the set of nodes of the crystal graph, and the graph itself, but will always make clear by context when we put the standard crystal structure on this set. We let $[n]^{\otimes k} $ denote the tensor product of $k$ copies of the standard $GL_n$ crystal with itself. Thus, nodes in the crystal $[n]^{\otimes k} $ are identified with words of length $k$ from the set $[n]$. The notation $[n]^{\otimes k} $ is used for both this set of words, and the crystal graphs whose nodes are labeled with the words. That is, we identify the node in the crystal $[n]^{\otimes k} $ given by the tensor product
\[ \omega_{1}\otimes \omega_{2}\otimes  \cdots \otimes \omega_{k-1}\otimes   \omega_k \in [n]^{\otimes k} \]
with the word
\[ \omega_{1}\omega_{2} \cdots \omega_{k-1} \omega_k. \]

There are two choices for defining crystal operators $e_i, f_i$ on $[n]^{\otimes k}$. We will adopt the \emph{Kashiwara} convention for crystal operators~\cite{Kash-Nak} (precise definitions in Section~\ref{[n] def} below). This is (essentially) the reverse to what is found in many sources, including Bump-Schilling~\cite{BumpSchilling}, but it has been used for some recent constructions (See~\cite{GerberLecouvey}). We do this because it allows us, when reading a word from left to right, to construct an associated perforated tableau, also obtained left to right.

\medskip

While $[n]^{\otimes k}$ decomposes as a sum of irreducible crystal graphs, and all isomorphism types of $GL_n$ crystal structures are found among them, other combinatorial models for crystal graphs such as \emph{biwords}, \emph{matrix models,} and \emph{RSK $(P,Q)$ pairs}, have been studied as well~\cite{BumpSchilling,FuLascoux,Fulton,GerberLecouvey,HeoKwon,Kash-Nak,Lee,NakayashikiYamada,Shimozono,vanLeeuwen}.

For example, start with some $\omega \in [3]^{\otimes 11}$:
\[ \omega=2122331331, \]
on which crystal operators $e_i$, $f_i$ act within a $GL_3$ crystal structure. We \emph{parse} $\omega$ into weakly decreasing factors, and denote the parsing by $\cal P$ and the parsed word as $\omega_{\cal P}$:
\[ \omega_{\cal P} = 21|22|331|331. \]

In Section~\ref{word to ptab} we define a bijection from parsed words to ptableaux. For now, we associate to $\omega_{\cal P}$ a \emph{biword} by putting constant factors composed of $1$'s, and then $2$'s, etc., over each successive factor of $\omega_{\cal P}$:

\[ \left( \begin{matrix} \tau \\ \omega \end{matrix} \right) = \left(\begin{matrix} 1 & 1& 2 & 2& 3 & 3& 3& 4&4&4 \\
2 & 1 &2 & 2& 3&3&1&3&3&1 \end{matrix}\right). \]

The $GL_3$ crystal structure on this biword is obtained simply by letting the crystal operators act on $\omega$. It can be shown that such operators preserve the parsing~\cite{BumpSchilling, Fulton}.

\medskip

In addition to biwords, \emph{matrix models} for crystal graphs have been studied~\cite{BumpSchilling,Fulton,NakayashikiYamada,Shimozono,vanLeeuwen}. Continuing our example, we compute a $4 \times 3$ \emph{matrix} $M$ associated to our biword, by setting the $(i,j)$ entry of $M$ equal to the number of times $i$ appears over $j$ in the biword:

\[ M = \begin{bmatrix} 1 & 1 & 0 \\ 0 & 2 & 0 \\ 1 & 0 & 2 \\ 1 & 0 & 2 \end{bmatrix}. \]
This map is a bijection from biwords $(\tau,\omega)^{T}$, with $\tau \in [\ell]^{\otimes k}$ and $\omega \in [n]^{\otimes k}$ to $\text{Mat}_{\ell \times n}({\mathbb N})$, the $\ell \times n$ matrices with non-negative integer entries. In this case, a $GL_3$ crystal structure can be defined for matrices $M$ by computing along the \emph{columns} of $M$ (See~\cite{BumpSchilling,vanLeeuwen}), and one shows this structure commutes with the map from biwords to matrices, yielding isomorphic crystal structures.
\medskip

Lastly, we can also associate a pair $(P,Q)$ of \emph{semistandard Young tableaux} (SSYT) to this biword by means of \emph{insertion algorithms}~\cite{BumpSchilling,Fulton}. Again, there is more than one such choice of algorithm. Here we use \emph{column} insertion on the word $\omega$, going left to right, using $\tau$ to determine the recording tableau. In our example, we obtain:

\[ P = \begin{ytableau} 1 & 1 & 1 & 2 & 2 & 3 & 3 \\
2 & 3 \\ 3 \end{ytableau}, \quad Q = \begin{ytableau} 1 & 1 & 2 & 3 & 4 & 4 & 4 \\ 2 & 3 \\ 3 \end{ytableau}. \]

The $GL_3$ crystal structure on SSYT $(P,Q)$ of the same shape is that induced by the structure on the word $\omega$ (one shows that the shape of the associated tableaux, and indeed the tableau $Q$, are fixed by this action).

\medskip

It can be shown (using our choices for crystal operators on words, matrix conventions, and insertion schemes) that the maps between these three models (biwords, matrices, RSK $(P,Q)$ pairs) are actually \emph{isomorphisms} of their crystal graph structures.

\medskip
Where do perforated tableau fit among these models? There are simply defined bijections from biwords to ptableaux, too. Starting from the biword given above, we work from left to right. The entries in $\tau$ in the top row will become the content (entries) in an associated ptableau $T$. If $i$ in $\tau$ is above some $j$ in $\omega$ of the biword, we put an $i$ in row $j$ of $T$, requiring such entries form a horizontal strip (and potentially inserting blanks, or ``perforations") so that the $T$ is column strict:

\[ \left( \begin{matrix} \tau \\ \omega \end{matrix} \right) = \left(\begin{matrix} 1 & 1& 2 & 2& 3 & 3& 3& 4&4&4 \\
2 & 1 &2 & 2& 3&3&1&3&3&1 \end{matrix}\right) \Rightarrow T=\begin{ytableau} \ & & 1 && 3 & 4\\ & 1 &2 & 2&& \\3 & 3 & 4 & 4 && \end{ytableau}. \]
It is not hard to show this map is a bijection. In Section~\ref{[n] def} we define crystal operators $e_i, f_i$ on perforated tableaux (``ptableaux"), and prove that these operators commute with those given by the crystal graph structure on $[n]^{\otimes k}$, and hence determining, (in general), an $GL_n$ crystal structure on ptableaux isomorphic to the structures of the other examples.

\medskip

   We hope to demonstrate (1) perforated tableaux are often among the \emph{simpler} models to use for various computations, including computations of crystal operators and determining highest weights, and (2) many combinatorial algorithms that have appeared in several \emph{distinct} models all naturally appear \emph{within} the perforated tableaux model.

   \medskip

Among the key advantages realized in perforated tableaux, we have:

\begin{enumerate}
\item Computing crystal operators on ptableaux is done with an easy, single, \emph{jeu de taquin} move, in which content moves down under the lowering operators $f_i$, and upwards with the raising operators $e_i$.
\item A ptableau is highest weight when its content is partition-shaped (unperforated). Thus, we can determine highest weights \emph{visually}, without computation, unlike the $[n]^{\otimes k}$, biwords, or matrix models.

\item The perforated tableaux model has a natural and simple tensor product operation.

\item We can give a simple combinatorial enumeration of the \LR\ rule for tensor product multiplicities by enumerating certain ptableaux. This rule directly generalizes classical ``\LR\ fillings.'' The ``word condition" of the classical case is not necessary generally, but it does locate one particular crystal graph within a fixed isomorphism class of irreducible crystal graphs.

\item A version of Sch\"{u}tzenberger evacuation is used to compute the Lusztig involution on the crystal structure of RSK $(P,Q)$ pairs. We prove this evacuation is actually a product of ptableaux crystal operators, so that evacuation on SSYT is viewed, in this context, as the operator taking a highest weight ptableau to its corresponding lowest weight.

\item Let $B_{\mu}$ and $B_{\nu}$ be irreducible crystals of highest weights $\mu$ and $\nu$, respectively. Lenart defined a non-combinatorial crystal isomorphism $\phi: B_{\mu} \otimes B_{\nu} \rightarrow B_{\nu} \otimes B_{\mu}$, called a commutator~\cite{lenartI,lenartII}. We prove a generalized form of tableau switching for ptableaux that computes $\phi$ combinatorially for arbitrary highest weight elements in tensor products of irreducible ptableaux graphs. Our algorithm generalizes related work of James and Kerber~\cite{J-K}, and also Benkart, Sottile, and Stroomer~\cite{Sottile}.)

\end{enumerate}
We should note that one key feature of the RSK $(P,Q)$ tableaux model is its use in solving the isomorphism and plactic equivalence problem~\cite{BumpSchilling}: if two nodes in a crystal graph are represented by pairs $(P,Q)$ and $(P', Q')$, respectively, then the nodes lie in isomorphic crystals precisely when $Q$ and $Q'$ are the same shape, and they lie in the \emph{same} irreducible crystal graph when $Q=Q'$.  They are plactically equivalent (in same location in isomorphic crystal graphs) when $P=P'$. In ~\cite{A-W RSK}, we present a beautiful solution to this problem, using ptableaux.

\subsection{Duality and Bicrystals}

Another important feature of the biword, matrix, and RSK pair models is the existence of a \emph{duality} map, allowing one to put a \emph{bicrystal} structure on these combinatorial objects (as $GL_{\ell} \times GL_{n}$ crystals for appropriate $\ell$), along with several fundamental results relating crystal operators on the $GL_{\ell}$ structure with combinatorially defined operations on the $GL_n$ crystal~\cite{BumpSchilling,FuLascoux,Fulton,GerberLecouvey,
NakayashikiYamada,vanLeeuwen}. Defining the dual to a perforated tableau is \emph{especially} simple, so we mention a few of these connections, taking up the general issue is outside the scope of this paper and will have to be taken up in later work.
\medskip

\begin{df} Let ${\mathcal Bi}([\ell]^{\otimes k},[n]^{\otimes k})$ denote the set of length $k$ \emph{biwords} over $[\ell]$ and $[n]$:

\[ B \in{\mathcal Bi}([\ell]^{\otimes k},[n]^{\otimes k}), \quad B = \left( \begin{matrix} a_1 & a_2 & \cdots & a_k \\ b_1 & b_2 & \cdots & b_k \end{matrix} \right) = \left( \begin{matrix} \tau \\ \omega \end{matrix} \right), \]
where $\tau \in [\ell]^{\otimes k}$ and $\omega \in [n]^{\otimes k}$. We restrict the words in the top row $\tau \in [\ell]^{\otimes k}$ to those that are \emph{weakly increasing} (so $a_i \leq a_{i+1}$), and we restrict words $\omega \in [n]^{\otimes k}$ to those that satisfy $b_i \geq b_{i+1}$ whenever $a_i = a_{i+1}$.
\end{df}

We say the \emph{dual} of a biword $(\tau,\omega)^{T} \in {\mathcal Bi}([\ell]^{\otimes k},[n]^{\otimes k})$ is the biword $(\widehat{\omega},\widehat{\tau})^{T}$, obtained by swapping the top and bottom rows of $(\tau,\omega)^{T}$, and reordering columns so that $\widehat{\omega}$ is weakly increasing, and under each constant factor of $\widehat{\omega}$ the corresponding factor of $\widehat{\tau}$ is weakly decreasing (that is, so that $(\widehat{\omega},\widehat{\tau})^{T} \in {\mathcal Bi}([n]^{\otimes k},[\ell]^{\otimes k})$. Thus, if
\[ \left( \begin{matrix} \tau \\ \omega \end{matrix} \right) = \left(\begin{matrix} 1 & 1& 2 & 2& 3 & 3& 3& 4&4&4 \\
2 & 1 &2 & 2& 3&3&1&3&3&1 \end{matrix}\right) \in {\mathcal Bi}([4]^{\otimes 11},[3]^{\otimes 11}), \]
we swap rows and reorder columns so that
\[ \left( \begin{matrix} \widehat{\omega} \\ \widehat{\tau} \end{matrix} \right) = \left( \begin{matrix} 1 & 1 & 1& 2 & 2&2&3&3&3&3\\ 4&3&1&2&2&1&4&4&3&3\end{matrix} \right) \in  {\mathcal Bi}([3]^{\otimes 11},[4]^{\otimes 11}). \]

A $GL_{3} \times GL_4$ \emph{bicrystal} structure is obtained on the pair $(\widehat{\omega}, \widehat{\tau})$ letting the $GL_3$ crystal operators act on $\omega$ and $GL_4$ operators act on $\widehat{\tau}$. One shows that these actions commute on the dual pairs $(\widehat{\omega}, \widehat{\tau})$.
\medskip

In the matrix setting, duality is even simpler to describe. We already had a $GL_3$ crystal structure defined by operators acting on the columns of the  $4 \times 3$ matrix $M$ given above. A $GL_{3} \times GL_4$ \emph{bicrystal} structure is obtained by defining the corresponding $GL_4$ structure applied to the columns of $M^{T}$, the transpose of $M$ (or appropriately defined versions of row operations on $M$). Again, one shows that these actions commute~\cite{BumpSchilling,vanLeeuwen}.

\medskip
For RSK $(P,Q)$ pairs, the dual of $(P,Q)$ is simply $(Q,P)$. As illustrated in our example
\[ P = \begin{ytableau} 1 & 1 & 1 & 2 & 2 & 3 & 3 \\
2 & 3 \\ 3 \end{ytableau}, \quad Q = \begin{ytableau} 1 & 1 & 2 & 3 & 4 & 4 & 4 \\ 2 & 3 \\ 3 \end{ytableau}, \]
$Q$ has content in $[4]$, and the induced $GL_4$ structure on $(Q,P)$ is given by acting on the word obtained from $Q$ by reading its rows from right to left, top to bottom.  One proves that this action commutes with the $GL_3$ crystal structure on $P$ described above.

\bigskip
We now determine the dual of a perforated tableau. Given a ptableau $T$,
its \emph{dual} $\widehat{T}$ by letting the content in row $i$ of $T$, read right to left, name the \emph{rows} in which a horizontal strip of $i$'s appears in $\widehat{T}$. Thus, if
\[ T=\begin{ytableau} \ & & 1 && 3 & 4\\ & 1 &2 & 2&& \\3 & 3 & 4 & 4 && \end{ytableau}\, , \]
the first row (read right to left) says that $\widehat{T}$ has a strip of $1$'s (moving left to right) appearing in rows $4$, $3$ and $1$, and so on:
\[ \widehat{T} = \begin{ytableau} \ & & 1 & &2 \\ &&2&2&\\ &1&&3&3\\1&3&3&& \end{ytableau}\, . \]
We obtain a $GL_4 \times GL_3$ bicrystal structure on the pair $(T, \widehat{T})$ by the induced $GL_3$ structure on $T$, and the $GL_4$ crystal structure on $\widehat{T}$.

Assuming facts stated for the other models above, we can prove the bicrystal structures on ptableaux commute. This follows from the easily proved bijection $\Phi$ from ptableaux to matrices: given $T$, the $(i,j)$ entry of the associated matrix $M$ is the number of times an $i$ appears in row $j$ of $T$.  For our example,
\[ T= \begin{ytableau} \ & & 1 && 3 & 4\\ & 1 &2 & 2&& \\3 & 3 & 4 & 4 && \end{ytableau}\stackrel{\Phi}{\Longrightarrow} M =  \begin{bmatrix} 1 & 1 & 0 \\ 0 & 2 & 0 \\ 1 & 0 & 2 \\ 1 & 0 & 2 \end{bmatrix}. \]
It is not hard to show that if, under this bijection, the ptableau $T$ corresponds to matrix $M$, then the ptableau $\widehat T$ corresponds to $M^T$, the transpose of $M$.

\[ \begin{array}{ccc}
T=\begin{ytableau} \ & & 1 && 3 & 4\\ & 1 &2 & 2&& \\3 & 3 & 4 & 4 && \end{ytableau}& \Longrightarrow& \widehat{T} = \begin{ytableau} \ & & 1 & &2 \\ &&2&2&\\ &1&&3&3\\1&3&3&& \end{ytableau}\\
 & & \\
\Phi \downarrow & & \Phi \downarrow \\
& & \\
M =  \begin{bmatrix} 1 & 1 & 0 \\ 0 & 2 & 0 \\ 1 & 0 & 2 \\ 1 & 0 & 2 \end{bmatrix} & \Longrightarrow &
M^{T}= \begin{bmatrix} 1&0 &1&1 \\1 & 2 & 0&0\\ 0&0&2 &2 \end{bmatrix}. \end{array} \]

Indeed, we can recover $\widehat{T}$ even more directly from the parsed word $\omega_{\cal}= 21|22|331|331$ that generated all our examples. The content of row $i$ of $\widehat{T}$ is merely the $i$-th factor of $\omega_{\cal P}$, read right to left. This content is then arranged minimally to ensure the resulting tableau filling is column strict.

A few comments are in order here. First, these observations emphasize that a perforated tableau $T$ associated to parsed word $\omega_{\cal P}$ is \emph{not} merely a rearrangement of the factors of $\omega_{\cal P}$, although its dual $\widehat{T}$, is. This is significant for a number of reasons. The number of rows of $\widehat{T}$ depends on the number of factors into which the associated word $\omega$ is parsed. That is, there is an \emph{intrinsic} $GL_n$ crystal structure on words $\omega \in [n]^{\otimes k}$, but the $GL_\ell \times GL_{n}$ bicrystal structure depends on the choice of parsing in the case of $\omega_{\cal P}$ or, equivalently, on the choice of compatible top row $\tau$ in the associated biword. One obtains different bicrystal structures, though all of them yield \emph{isomorphic} $GL_n$ structures.

We will see that one consequence of this is that if $T$ is the ptableau associated to a parsed word $\omega_{\cal P}$, and $T'$ is the ptableau associated to $\omega_{\cal P'}$ (the same word, but with a different parsing), then $T$ and $T'$ will be \emph{plactically equivalent}, meaning that they lie in the same location in isomorphic irreducible crystal graphs, and (by means of our Theorem~\ref{highest weight} below) their respective highest weight elements will be partitions of the \emph{same shape}. This need not be true of their duals $\widehat{T}$ and $\widehat{T'}$, which might not even have the same number of rows, despite both being more directly determined by the word $\omega$.

\medskip
We should mention that in recently available work of Gerber and Lecouvey~\cite{GerberLecouvey} (which appeared subsequent to the authors' earlier draft of this work~\cite{A-W Ptab}), a crystal model consisting of a tensor product of column tableaux is proposed.  They define a crystal operator from a word in $[n]^{\otimes k}$ computed from their model, and give many other significant results, including a definition of duality that allows them to compute the maps mentioned above.  Most of these results are outside the scope of our work here.

\section{Definitions, Perforated Tableaux}

In all that follows our Young diagrams (or just diagrams) will be rectangular. A \emph{tableau} is a filling of boxes of the diagram with positive integers, possibly leaving some boxes blank (or \emph{perforated}). We refer to boxes in a tableau that are not blank as \emph{having content}.

\begin{df} Given a box $\framebox{b}$ in a tableau $T$, the \emph{northwest shadow} of $\framebox{b}$ is the collection of boxes in $T$ to the west, north, or both, of $\framebox{b}$.
\end{df}

\begin{df} In a tableau, a collection of boxes forms a \emph{horizontal strip} if no two boxes in the strip are in the same column, and any box in row $k$ of the strip is in a column to the right of any box of the strip in row $k+\ell$, for $\ell > 0$. Given a horizontal strip, we call the lowest, left-most box the \emph{head} of the strip, and the highest, rightmost box the \emph{tail}. If the collection of all the appearances of an integer $i$ in a tableau form a horizontal strip, we call this the \emph{$i$-strip} of the tableau.
\end{df}

In this paper, we picture horizontal strips as if they ``grow'' from the lower left to right, bottom to top (from head to tail). This will seem natural when we introduce a map from words in $[n]^{\otimes k}$ to perforated tableaux, in which strips are filled from head to tail as the word is read in (via our map) from left to right.

\ytableausetup{boxsize=1em}
\begin{df} \label{ptab} A rectangular tableau $T$ will be called a \emph{perforated tableau}  (or \emph{ptableau}, for short) if:
 \begin{enumerate}
 \item All boxes in $T$ contain positive integers (have content), or are blank (unfilled, denoted $\Box$).
  \item For any positive integer $i$, the $i$'s in $T$ form a horizontal strip.
 \item For positive integers $i$ and $j,$ with $i < j$, every element in the horizontal strip of $j$'s lies \emph{outside} the northwest shadow of any entry in the horizontal strip of $i$'s.
 \item $T$ has no column entirely composed of blanks (rows of all blanks are allowed).
\end{enumerate}

\noindent If, in addition to the above requirements, a perforated tableaux $T$ satisfies:

\medskip

5. For any positive integer $i$ appearing in $T$, the head of the $i$-strip of $T$  lies in a row \emph{strictly} below the row containing the tail of $(i-1)$-strip, then we say $T$ is a \emph{minimally parsed} perforated tableau.
\end{df}

Here is an example of a minimally parsed perforated tableau:

\ytableausetup{smalltableaux}
 \ytableausetup{centertableaux}
\[T = \begin{ytableau}
\ &&&&&&&&1&1&1\\
&&&&&1&1&1&&2&2\\
&&&1&1&&2&2&3&3&3\\
1&1&1&2&2&2&3&3&4&4& \end{ytableau}, \]
with $\wt(T) = (3,5,7,10)$.

\medskip
\noindent \emph{Notational Convention:} We need to distinquish between an actual box of a ptableau, and the \emph{value} (an integer or blank) of that box. We denote a box with content $a$ by $\framebox{a}$; we denote a blank (empty) box by $\Box$.

\begin{df}
We define the \emph{weight} of a perforated tableau $T$ to be the sequence $\wt(T) = (m_{1}, \ldots , m_{n})$, where $m_{i}$ is the number of filled boxes in row $i$ of $T$.
\end{df}

\begin{df} Let $\nu=(\nu_1, \ldots , \nu_n)$ be a partition (so $\nu_{1} \geq \nu_{2} \geq \cdots \geq \nu_{n}$). Let $\SSYT(\nu)$ denote the set of semistandard Young tableaux of shape $\nu$. All tableaux in $\SSYT(\nu)$ will have $\nu_{1}$ columns, with content in every box of row $1$. In row $i$, for $2 \leq i \leq n$, there is content in the first $\nu_{i}$ boxes, followed by $(\nu_{1} - \nu_{i})$ blanks.
\end{df}

So, an example of a semistandard Young tableau of weight $(4,3,3,1)$ is:
\[\begin{ytableau} 1 & 1 & 2 & 2 \\ 2 & 2 &3 & \\ 3 & 4 & 4 & \\4 & & & \end{ytableau} , \quad
\hbox{while one of weight $(4,3,3,1,0)$ is} \quad
\begin{ytableau} 1 & 1 & 2 & 2 \\ 2 & 2 &3 & \\ 3 & 4 & 4 & \\4 & & & \\ &&& \end{ytableau} . \]
Note that, typically, a SSYT is drawn without including the blank boxes, but, by including the blank boxes, we see that SSYT are a special case of perforated tableaux.

\begin{df} Suppose $S$ and $T$ are tableaux with $n$ rows. We say $S$ and $T$ are \emph{row equivalent} if we can transform $S$ into $T$ by a finite number of steps in which a blank in a given row is swapped with an adjacent non-blank in the same row, such that the result is a valid perforated tableau. Note that if $S$ and $T$ are row equivalent, then $\wt(S) = \wt(T)$.
\end{df}

Row equivalence is clearly an equivalence relation on tableaux.

\medskip

\begin{df}
A ptableau $T$ is \emph{left-justified} if all of $T$'s content is as far to the \emph{left} as possible. That is, in any ptableau $T'$ that is row equivalent to $T$, each entry of $T$ is weakly to the left of the corresponding entry in $T'$.

Similarly, a ptableau $T$ is \emph{right-justified} if all content of $T$ is as far right as possible.
\label{justified}
\end{df}

The following are then easily proved:

\begin{lem}
\begin{enumerate}
\item Given any ptableau $T$, there is a unique ptableau $^*T$ that is left-justified and row equivalent to $T$.
\item Given a left-justified ptableau $^*T$, there is a unique right-justified ptableau $T^*$ that is row equivalent to $^*T$.
\end{enumerate}\label{left-right-justified}
\end{lem}

\begin{df} Let $\PTab$ denote the set of row equivalence classes of ptableaux, each of which contains ptableaux that are row-equivalent to the representative $^*T$ and that have the same number of columns as $^*T$.  Let $\PTab_{n}$ denote the subset of $\PTab$ composed of ptableaux with $n$ rows (some of those rows may be entirely blank).

Let $\PTab_{(\ell,n)}$ denote the subset of $\PTab_{n}$ whose content is from the set $[\ell]$.
\label{ptab def}
\end{df}

We require that all members of an equivalence class have the same number of columns to rule out such examples as:

\[ T=\begin{ytableau} \ &&& & 1& 1& 1 & 2 \\ 1&2&2&3&&& & \end{ytableau}\, . \]

Strictly speaking, this is a perforated tableau, but it would seemingly be row-equivalent to:
\[ \begin{ytableau} \ & &1& 1& 1 & 2&& \\ 1&2&2&3&&& & \end{ytableau}\, , \]
were we to preserve the blanks in $T$.  But this is not a perforated tableau because by definition, ptableaux cannot have columns of blanks.  We resolve this by using the left-justified form of $T$, with blank columns removed:

\[ ^*T = \begin{ytableau} \ & 1& 1& 1 & 2 \\ 1&2&2&3&  \end{ytableau} \]
and then only admit into the equivalence class of $^*T$ those perforated tableaux with the same number of columns (in this case, 5,) as $^*T$.

Note that, by Lemma~\ref{left-right-justified}, the right-justified ptableau
\[ T^* = \begin{ytableau} \ & 1& 1& 1 & 2 \\ 1&&2&2&3 \end{ytableau} \]
is in the same \PTab\ equivalence class as $^*T$, but that this class would \emph{not} include the tableau:

\[  \begin{ytableau} \ & 1& 1& 1 & 2& \\ 1&2&2&&&3  \end{ytableau} \]
because it has more columns than $^*T$.

We regard row-equivalence classes of ptableaux, under Definition~\ref{ptab def}, as a single object, which (to avoid unnecessary circumlocution) we will still refer to as ``a'' ptableau. We think of this as a single tiled grid (with the number of columns defined by $^*T$) in which the tiles (the content) are free to slide right or left, subject to the constraint of column strictness. Thus, we regard

\[\begin{ytableau}
\ &&&&1&1&4 \\
&&1&1&&2&5\\
1&2&3&4&4&5& 6\end{ytableau}\ \  \text{ and } \ \ \begin{ytableau}
\ &&&1&1&4 &\\
&1&1&2&&&5\\
1&2&3&4&4&5&6 \end{ytableau} \]
to be equal as ptableaux since they are row equivalent.  Furthermore,  \[  \begin{ytableau}
\ &&&1&1&4 &\\
&1&1&2&&&5\\
1&2&3&4&4&5&6 \end{ytableau} \]
is the (unique) left-justified representative of this equivalence class. So, while a ptableau $T$ can potentially be represented in multiple (row-equivalent) forms, $^*T$ and $T^*$ are unique representatives.

\section{From Words $[n]^{\otimes k}$ and Biwords to Perforated Tableaux}\label{word to ptab}

\begin{df} Given a word $\omega \in [n]^{\otimes k}$, a \emph {subword} of $\omega$ is sequence of letters in $\omega$, in order but not necessarily contiguous.  A \emph {factor} of $\omega$ is a contiguous subword of $\omega$.
\end{df}

Below we define a map $\Pf$ from words $[n]^{\otimes k}$ to perforated tableaux. The map will be determined by a word $\omega \in [n]^{\otimes k}$, and a ``parsing'' of $\omega$ into weakly decreasing factors. Condition (5) of Definition~\ref{ptab} above characterizes the form of the ptableau $\Pf(\omega)$ when $\omega$ has been parsed into the \emph{minimal} number of such factors, though condition (5) can be verified independently.

\begin{df} Given some word
$\omega = \omega_{1}\omega_{2} \cdots \omega_{k-1} \omega_k \in [n]^{\otimes k} ,$
a \emph{parsing} of the word is any partition of $\omega$ into a disjoint collection of weakly decreasing factors, when read left to right. Trivial parsings (with empty factors) are allowed.
The \emph{minimal parsing} of $\omega$ is the (unique) parsing of $\omega$ into the fewest possible number of weakly decreasing factors.
\end{df}
For example, if $\omega$ is
\[ 4433442233344411122233444 \in [n]^{\otimes 24}, \]
we parse it into weakly decreasing factors by denoting the breaks between factors with the symbol ``$\, |$":
\[ 4433|4422|333|44|4111|2222|33||444. \]
(the ``$\, ||$" above denotes an empty factor).  The minimal parsing is:
\[  4433|4422|333|444111|2222|33|444. \]

For future reference, we record the following:

\begin{df}
Given some word $\omega \in  [n]^{\otimes k} $, we denote a parsing of $\omega$ by the symbol ${\cal P}$, and denote the word $\omega$, with its associated parsing, by $\omega_{\cal P}$. The minimal parsing of a word $\omega$ is denoted $\omega_{{\cal P}_{\mmin}}$.

We let $P([n]^{\otimes k})$ denote the set of all parsings of words in $[n]^{\otimes k}$.

If a parsed word $\omega_{\cal P} \in P([n]^{\otimes k})$, breaks $\omega$ into $\ell$-many weakly decreasing factors, let ${\mathcal Bi}(\omega_{\cal P})$ denote the unique biword in ${\mathcal Bi}([\ell]^{\otimes k},[n]^{\otimes k})$ obtained by writing a string of $i$'s over the $i$-th factor.\label{parse def}
\end{df}

 Thus, our example
\setcounter{MaxMatrixCols}{30}

\[ \omega_{\cal P} = 4433|4422|333|44|4111|2222|33||444 \]

would be associated to the biword

\[ {\mathcal Bi}(\omega_{\cal P}) = \left( \begin{matrix} 1&1&1&1&2&2&2&2&3&3&3&4&4&5&5&5&5&6&6&6&6&7&7&9&9&9\\
4&4&3&3&4&4&2&2&3&3&3&4&4&4&1&1&1&2&2&2&2&3&3&4&4&4 \end{matrix} \right). \]

Consequently, we have the remarkable fact that parsed words carry the \emph{same} information (in the context of models for crystal operators) as biwords, RSK $(P,Q)$ pairs, and matrix models $\text{Mat}_{\ell \times n}$. Below, we define a map from parsed words to ptableaux. It will be obvious that we could have defined them from biwords directly, but it will be both notationally convenient to use parsed words, and it will allow us to state with precision how the crystal structure we define on ptableaux commutes with this map, using the crystal operators on words (and ignoring the parsing).

\medskip

We define the map $\Pf:P([n]^{\otimes k}) \rightarrow \PTab_{n}$ taking a parsed word $\omega_{\cal P} \in P([n]^{\otimes k})$ to an equivalence class $\Pf(\omega_{\cal P}) = T \in\PTab_{n}$ determined by its unique left-justified perforated tableau $^*T$, constructed from $\omega$ and its parsing ${\cal P}$. It will be apparent that $^*T$ clearly fulfills the requirements of Definition~\ref{ptab}, so we omit the proof of this claim.

Given a word
$ \omega_{\cal P} = \omega_{1}\omega_{2} \cdots \omega_{k-1} \omega_k \in P([n]^{\otimes k}) ,$
we represent the parsing of $\omega$ by:

\[ \omega_{\cal P}= h_{1} | h_{2}|  \cdots|  h_\ell, \]
where each $h_{s}$ is a weakly decreasing factor of $\omega_{\cal P}$. Given a factor $h_{s}$, let $h_{s}^{i}$ denote the number of $i$'s in $h_{s}$.

We construct from $\omega_{\cal P}$ a left-justified ptableau $^*T$:  Begin with a blank tableau of $n$ rows (and an as-yet unspecified number of columns), with the left-most column designated as column 1, the next column as column 2, etc.   Given the factor $h_{1}$, we put $h_{1}^{n}$-many $1$'s in row $n$, starting in column 1, then $h_{1}^{n-1}$-many $1$'s in row $n-1$, starting in column $h_{1}^{n}+1$, then $h_{1}^{n-2}$-many $1$'s in row $n-2$, starting in column $h_{1}^{n} +h_{1}^{n-1} +1$, etc., until we have constructed a horizontal strip of $1$'s. We then put $h_{2}^{k}$-many $2$'s in row $k$, starting from row $k=n$ and moving up, subject to the constraint that each row of $2$'s is as far left as possible,  and the filling is column-strict (leaving blanks in the tableau as needed to achieve this). We continue in the same manner for $h_{3}$ through $h_{\ell}$. The rightmost column containing an integer is then the rightmost column of the perforated tableau.
We denote the output of this construction by $^*T$ since it is clearly left-justified. The ptableau $^*T$ determines, by definition, a row-equivalence class $\Pf(w_{\cal P}) \in \PTab$.

\medskip

As an example, let

\[ \omega= 44433222111444333221443344 \in [4]^{26}.\]
The minimal parsing of $\omega$ is

\[ \omega_{{\cal P}_{\mmin}}= \underbrace{44433222111}_{h_1}|\underbrace{444333221}_{h_2}|
\underbrace{4433}_{h_{3}}|\underbrace{44}_{h_{4}}. \]


\noindent  To construct $^*T \in \Pf(\omega_{{\cal P}_{\mmin}})$, we first construct the horizontal strip of $1$'s determined by $h_1=44433222111$. Working left to right in $h_{1}$, we place three $1$'s in row $4$, then two $1$'s in row 3, three $1$'s in row 2, and finally three $1$'s in row 1:

\[ \begin{ytableau}
\ &&&&&&&&1&1&1\\
&&&&&1&1&1&&&\\
&&&1&1&&&&&&\\
1&1&1 &&&&&&&&\end{ytableau}. \]

\noindent We then insert a maximally left-justified horizontal strip of $2$'s determined by $h_2$. Notice that this increases the number of columns in the ptableau.
\[ \begin{ytableau}
\ &&&&&&&&1&1&1&2\\
&&&&&1&1&1&&2&2&\\
&&&1&1&&2&2&2&&&\\
1&1&1&2&2&2& &&&&&\end{ytableau}. \]

\noindent We fill in the rest of the perforated tableau using $h_3$ and then $h_4$, giving the left-justified $^*T \in \Pf(\omega_{{\cal P}_{\mmin}})$ where:

\[ ^*T=\begin{ytableau}
\ &&&&&&&&1&1&1&2\\
&&&&&1&1&1&&2&2&\\
&&&1&1&&2&2&2&3&3&\\
1&1&1&2&2&2&3&3&4&4&& \end{ytableau}. \]

 In terms of a biword (induced by some parsed word $\omega_{\cal P})$:
\[ \left( \begin{matrix} a_1 & a_2 & \cdots & a_k \\ \omega_1 & \omega_2 & \cdots &\omega_k \end{matrix} \right), \]
we obtain $\Pf(\omega_{\cal P})$ from this biword by placing each entry $a_{i}$ in the top row of the biword into row $\omega_{i}$ of the ptableau, moving content to the right, as necessary, to avoid column strictness violations.

\begin{df} Define ${\mathcal Bi}\Pf: {\mathcal Bi}([\ell]^{\otimes k},[n]^{\otimes k}) \rightarrow \PTab_{n}$
via
\[ {\mathcal Bi}\Pf \left( \begin{matrix} \tau \\ \omega \end{matrix} \right) = \Pf (\omega_{\cal P}),
\]
where
\[ \omega_{\cal P}= {\mathcal Bi}^{-1}\left( \begin{matrix} \tau \\ \omega \end{matrix} \right) = \Pf (\omega_{\cal P}),
\]
and ${\mathcal Bi}$ is defined as in Definition~\ref{parse def}.
\end{df}



\begin{lem} Let $\Pf$ and ${\mathcal Bi}\Pf$ be the maps defined above.  Then recalling Definition~\ref{ptab def}:
\begin{enumerate}
\vskip .1in
\item ${\mathcal Bi}\Pf$ is a bijection from biwords $\bigcup_{k,\ell \geq 1}{\mathcal Bi}([\ell]^{\otimes k},[n]^{\otimes k})$ to equivalence classes in $\PTab_{n}$, where, in particular, ${\mathcal Bi}\Pf$ maps elements in ${\mathcal Bi}([\ell]^{\otimes k},[n]^{\otimes k})$ to $\PTab_{(\ell,n)}$.
\item Regarding each element $\omega \in [n]^{\otimes k}$ as a minimally parsed element of $P([n]^{\otimes k})$, then $\Pf$ is a bijection from $\bigcup_{k\geq 1}[n]^{\otimes k}$ to row equivalence classes of minimally parsed ptableaux in $\PTab_{n}$.
\item  The image $\omega_{\cal P} \mapsto  T=\Pf(\omega_{\cal P} )$ determines a bijection from letters in $\omega$ to the content of the ptableau $T$, because row equivalence fixes the content and order of entries in any row. We also denote this bijection as ``$\Pf$".
\item If a letter $\omega_{i}$ lies to the right of a letter $\omega_{j}$ in a word $\omega_{\cal P} \in [n]^{\otimes k}$, then $\Pf(\omega_{j})$ lies outside the northwest shadow of $\Pf(\omega_{i})$ in $T$.
\end{enumerate} \label{Perf order}
\end{lem}
Given the bijections between biwords, elements of $\text{Mat}_{\ell \times n}$ and RSK $(P,Q)$ pairs, we could extend bijections between ptableaux and these objects as well.

We will need the following technical result in the next section:

\begin{lem} Let $T \in \Pf(\omega_{\cal P})$ for some $\omega_{\cal P} \in P([n]^{\otimes k})$. Then the number of columns in $^*T$ equals the length of the longest weakly decreasing subword in $\omega_{\cal P}$.
\label{width lemma}
\end{lem}
\begin{proof} Let $m$ equal the length of the maximal weakly decreasing subword of $\omega_{\cal P}$ (in this setting, the parsing is irrelevant).  Recall the bijection between the entries in $\omega_{\cal P}$ and the content of $^*T\in \Pf(\omega_{\cal P})$ given in Lemma~\ref{Perf order}.  Any weakly decreasing subword of $\omega_{\cal P}$ will correspond, under this bijection, to distinct boxes of content in $^*T$ appearing in different columns.  Hence the number of columns of $^*T$ is at least $m$. Clearly, the number of columns cannot exceed $m$; if it did,  there would be a sequence of entries in $^*T$ appearing in distinct columns that would correspond to a weakly decreasing subword in $\omega_{\cal P}$ of length greater than $m$, contradicting $m$'s maximality.
\end{proof}

\section{Crystal Structures on $P([n]^{\otimes k})$ and $\PTab_n$}\label{[n] def}
\subsection{Crytal Operators on $P([n]^{\otimes k})$}

We briefly review the crystal graph structure on $[n]^{\otimes k}$ (the parsings are irrelevant to these considerations), following the Kashiwara convention on crystal operators. See~\cite{BumpSchilling}  Lemma 2.33, p.\ 21, (defined there with an alternate orientation from ours):

\begin{prop}
Suppose $\omega = \omega_{1}\omega_{2} \cdots \omega_k \in [n]^{\otimes k}$. For any $s \in [n]$, define $\phi_{i}(s) = \delta_{i}(s)$ (Kronecker delta function), and $\epsilon_{i}(s) = \delta_{i+1}(s)$. Then define

\begin{align*} c_{e}^{i}(\omega,j)  &=   \left( \sum_{s=1}^{j}\epsilon_{i}(\omega_{s}) - \sum_{s =1}^{j-1}\phi_{i}(\omega_{s}) \right),\\
c_{f}^{i}(\omega,j)&=\left( \sum_{s=j}^{k} \phi_{i}(\omega_{s}) - \sum_{s = j+1}^{k} \epsilon_{i}(\omega_{s}) \right), \end{align*}
and, for $1 \leq i <n$,
\[ \epsilon_{i}(\omega) = \max_{j=1}^{k}c_{e}^{i}(\omega,j), \qquad \phi_{i}(\omega) = \max_{j=1}^{k} c_{f}^{i}(\omega,j). \]

Then for $s\in [n]$, define $e_{i}(s) = i$ if $s= i+1$, and $\NULL$ otherwise, and $f_{i}(s) = i+1$ if $s=i$, and $\NULL$ otherwise. Extend these definitions to the word $\omega \in [n]^{\otimes k}$ by setting
\[ e_{i}(\omega) = \omega_{1} \omega_{2}  \cdots \omega_{j-1} e_{i}(\omega_{j})\omega_{j+1} \cdots \omega_k, \]
with $j$ the \emph{smallest} index at which $\epsilon_{i}(\omega) = c_{e}^{i}(\omega,j)$, and
\[ f_{i}(\omega) = \omega_{1} \omega_{2} \cdots \omega_{j-1}  f_{i}(\omega_{j})\omega_{j+1} \cdots \omega_k, \]
with $j$ the \emph{largest} index at which $\phi_{i}(\omega) =  c_{f}^{i}(\omega,j)$, or $\NULL$ if $e_{i}(\omega_{j})$ or $f_{i}(\omega_{j})$ is $\NULL$.

These definitions determine an $A_{n-1}$ crystal structure on $[n]^{\otimes k}$, as a $k$-fold tensor product of the standard crystal $[n]$. \label{crystal def}
\end{prop}

\begin{df} A \emph{crystal graph} on $[n]^{\otimes k}$ is a subset $B \subseteq [n]^{\otimes k}$ that is closed under the crystal operators $e_{i}, f_{i}$, for $1 \leq i < n$. We regard the words in $B$ as the nodes of the graph; two words $w, \omega' \in B$ are connected by an edge if $\omega' = e_{i}(\omega)$ for some $i$ or, equivalently, if $\omega = f_{i}(\omega')$.
A \emph{connected} crystal graph is called \emph{irreducible}.
\end{df}
\begin{df} A word $\omega \in [n]^{\otimes k}$ is \emph{highest weight} if and only if $e_{i}(\omega) = \NULL$ for all $i$, $1 \leq i < n$. Likewise, a word $\omega' \in [n]^{\otimes k}$ is \emph{lowest weight} if and only if $f_{i}( \omega') = \NULL$ for all $i$, $1 \leq i < n$.
\end{df}

 It is known that if $\omega \in  [n]^{\otimes k}$ is highest weight, then its weight, $\wt(\omega)$, is a partition. (Recall that the {\em weight} of a word $\omega$ is the composition $(a_1, a_2, \dots, a_n)$, where $a_i = $ the number of $i$'s in $\omega$.) Every irreducible crystal graph possesses a unique element of highest weight (\cite{BumpSchilling}, p.\ 34). Thus we need to  distinguish between a \emph{highest weight element} of an irreducible crystal graph, and the actual \emph{weight} of the highest weight element. This leads to the following:

\begin{df} Suppose $\omega \in [n]^{\otimes k}$ is highest weight. We will denote the connected component containing $\omega$ (an irreducible crystal graph) by ${\mathbb B}^{\nu}(\omega)$, where $\nu=\wt(\omega)$. \label{word irred}
\end{df}

 If $\omega, \omega' \in [n]^{\otimes k}$ are highest weight words such that $\wt(\omega)=\wt(\omega') =\nu$, then the irreducible crystal graphs ${\mathbb B}^{\nu}(\omega)$ and ${\mathbb B}^{\nu}(\omega')$ are necessarily isomorphic (see~\cite{BumpSchilling}, p. 52). We refer to words identified by means of such an isomorphism as \emph{plactically equivalent}.

\subsection{Crystal Operators on Perforated Tableaux}
\begin{df}
Let $T \in \PTab_{n}$  and choose $i$, for $1 \leq i <n$. Let $T[i,i+1]$ denote the ptableau formed by the entries of $T$ in rows $i$ and $i+1$.\end{df}

Note that, as columns of blanks are not allowed in a ptableau, we omit any such columns in the determination of $T[i,i+1]$. For example, if
\[ T=  \begin{ytableau}
\ &&&&&&&&4&4&5\\
&&&&1&1&2&&&&6\\
&&1&1&2&&&4&5&6&7\\
1&1&2&3&3&3&4&6&6&& \end{ytableau}\, , \quad \hbox{then} \ \  T[2,3] = \begin{ytableau}
\ &&1&1&2&&6\\
1&1&2&4&5&6&7 \end{ytableau} \, . \]

\begin{lem} Let $T \in PTab_n$. If content in row $(i+1)$ of $^*T[i,i+1]$ is uncovered, then the corresponding content in $^*T$ is also uncovered. \label{blank safe}
\end{lem}

\begin{proof} Replace rows $i$ and $i+1$ in $^*T$ with
$^*T[i,i+1]$, with its content as far left as possible.
\[\hbox{For example:} \quad \begin{ytableau}
\ &&&&&&&&4&4&5\\
&&1&1&2&&6&&&&\\
1&1&2&4&5&6&7&&&&\\
1&1&2&3&3&3&4&6&6&& \end{ytableau} \, .\]

The resulting tableau will typically no longer be a valid ptableau, but we shall sequentially adjust content in it to resolve this in such a way that (1) the end result is $^*T$, and (2) content uncovered in $^*T[i,i+1]$ will remain uncovered during the adjustment. Note that we have not changed the \emph{content}, nor the \emph {order} of the non-blank entries in rows $i$ and $i+1$ of $^*T$, and have only changed their relative positions in the rows.

Begin with the $1$'s in row $i+1$. If the rightmost $\framebox{1}$ of $^*T$ appearing below row $i+1$ is in column $\ell$, then move all content of rows $i$ and $i+1$ to the right, so that the first entry in row $i+1$ is in column $\ell+1$.  In this example, $\ell =2$ and we move all content in rows $2$ and $3$ two columns to the right:

\[ \begin{ytableau}
\ &&&&&&&&4&4&5\\
&&1&1&2&&6&&&&\\
1&1&2&4&5&6&7&&&&\\
1&1&2&3&3&3&4&6&6&& \end{ytableau}\rightarrow
\begin{ytableau}
\ &&&&&&&&4&4&5\\
&&&&1&1&2&&6&&\\
&&1&1&2&4&5&6&7&&\\
1&1&2&3&3&3&4&6&6&& \end{ytableau}\]

 Call this the ``$1$-adjustment''. If, after this move, the rightmost $\framebox{2}$ appearing below row $i+1$ appears in column $\ell'$, then we move \emph{all} content in rows $i$ and $i+1$ with values greater than $1$ to the right until the resulting (partially adjusted) ptableau is column-strict and left-justified in columns $1$ through $\ell'$. Call this shift the ``$2$-adjustment''. (In our example, the 2-adjustment isn't necessary as the 2's in rows $3$ and $4$ are already column-strict.)  We continue inductively until we have adjusted all content from its original location in $^*T[i,i+1]$ to where it must appear in $^*T$, in order to recreate column-strictness, maintaining the requirement that if an $s$-adjustment moves content in rows $i$ and $i+1$, it moves all content with values greater than or equal to $s$ an equal number of columns to the right.  In our example, performing the $4$-adjustment moves all content greater than or equal to $4$ in rows $2$ and $3$ two columns to the right and results in $^*T$:
 \[ \begin{ytableau}
\ &&&&&&&&4&4&5\\
&&&&1&1&2&&6&&\\
&&1&1&2&4&5&6&7&&\\
1&1&2&3&3&3&4&6&6&& \end{ytableau}\text {   4-adjustment:}\rightarrow
\begin{ytableau}
\ &&&&&&&&4&4&5\\
&&&&1&1&2&&&&6\\
&&1&1&2&&&4&5&6&7\\
1&1&2&3&3&3&4&6&6&& \end{ytableau}.\]

Now let $\framebox{c}$ denote an uncovered entry in row $(i+1)$ of $^*T[i,i+1]$. All content in row $i$ in a column strictly left of the column containing $\framebox{c}$ has value smaller than $c$, and content in row $i$ to the right has value greater than or equal to $c$; otherwise, we could have covered $\framebox{c}$ in $^*T[i,i+1]$. However, by the process noted above, at each stage as we adjust $^*T$ from $^*T[i,i+1]$, for $s \leq c$, the $s$-adjustment moves all content with values greater than or equal to $s$ the same number of entries to the right. In particular, if $\framebox{c}$ is uncovered before such an $s$-adjustment, it remains uncovered after it. For any $s>c$, the $s$-adjustment will not move $\framebox{c}$ at all. Thus, no content in row $i$ will be moved over $\framebox{c}$ in $^*T$, and the lemma is proved.
\end{proof}

 \medskip


\begin{df} \label{e_i def} Suppose $T \in \PTab_{n}$ and that $\framebox{c}$ is an entry in row $i+1$ of $T$ such that in $^*T[i,i+1]$,  $\framebox{c}$ is the rightmost entry with a blank above it. In this case, we define $e_i (T)$ to be the ptableau obtained by swapping $\framebox{c}$ with the blank above it (in row $i$) of $^*T$ (by virtue of Lemma~\ref{blank safe}). If there is \emph{no} such entry $\framebox{c}$ in ${^*T[i,i+1]}$ , then $e_i (T) = NULL$.

\noindent We define $\epsilon_{i}(T)$ to be the number of blanks in row $i$ (the top row) of $^*T[i,i+1]$.
\end{df}

 We give an analogous definition of operators $f_{i}: \PTab_{n} \rightarrow \PTab_{n} \cup \{ \NULL \}$ below.

 \medskip

Recall the definition of the \emph{dual} $\widehat{T}$ of a ptableau $T$ given in the introduction. There are related results~\cite{Shimozono,GerberLecouvey} that use the bicrystal structure on a crystal graph model to realize a crystal operator on one object as some \emph{other} combinatorially defined transformation on its \emph{dual}.  This includes using \emph{jeu de taquin} operations on dual elements, to compute commutators of tensor products of factors of words.  That is not (quite) what is asserted here. Rather, in the case of crystal operators $e_i$ and $f_i$ on ptableaux, the operators \emph{themselves} can be seen as a sort of \emph{jeu de taquin} operation on the ptableau $T$. This may very well correspond to an interesting operation on the dual $\widehat{T}$. However, those previously cited results do not imply our theorems here, and second, our goal is to \emph{replace} the more complicated definitions of crystal operators (on words, matrices, SSYT) with the simpler operators defined on ptableaux.

Computing these (to be shown) crystal operators on ptableaux generally requires less computation than computing them on words. Computation of the crystal operators $f_{i}, e_{i}$ on words in $[n]^{\otimes k}$ requires counting $i$'s and $(i+1)$'s across the \emph{entire} length of the word (computing in the other models is not substantially simpler), while the computation of operators on ptableaux requires working only with the columns of rows $i$ and $i+1$ in $T$ that contain a blank in row $i$.

\subsection{Examples}

As an example, consider $T$ and $T[2,3]$ as above.

We calculate $e_{2}(T)$:
\begin{equation} ^*T[2,3] = \begin{ytableau}
\ &&1&1&2&&6\\
1&1&2&4&5&6&7 \end{ytableau} ,\label{example eq} \end{equation}
and so the $6$ in the sixth column of $^*T[2,3]$ is the rightmost entry under a blank. Noting Lemma~\ref{blank safe} (ensuring the blank persists in $^*T$), we circle this entry in $^*T$ below:
\ytableausetup{centertableaux}
\[ ^*T=  \begin{ytableau}
\ &&&&&&&&4&4&5\\
&&&&1&1&2&&&&6\\
&&1&1&2&&&4&5&\textcircled{$\scriptstyle 6$}&7\\
1&1&2&3&3&3&4&6&6&& \end{ytableau}\, ,
\quad \hbox{and so} \quad
e_{2}(T) =  \begin{ytableau}
\ &&&&&&&&4&4&5\\
&&&&1&1&2&&&\textcircled{$\scriptstyle 6$}&6\\
&&1&1&2&&&4&5&&7\\
1&1&2&3&3&3&4&6&6&& \end{ytableau}\, . \]
By definition, $\epsilon_2(T) = 3$ because there are three blanks in the top row of $T[2,3]$.

Thus, $T \in \Pf (\omega)$ for the word:
$ \omega_{{\cal P}}= 443322|432|444|4311|31|4432|3.$

Anticipating Theorem~\ref{ei map} below, we note:
$e_{2}(\Pf(\omega_{{\cal P}})) = \Pf(e_{2}\left(\omega_{{\cal P}}\right)). $
\medskip

The definitions of crystal operators on ptableaux make many features of the theory of crystal graphs more visually clear. For example, the form of the $1$-root string (see~\cite{BumpSchilling}, p.23),

\[322|3311|222|3\xrightarrow{e_1}322|3311|221|3 \xrightarrow{e_1}321|3311|221|3\xrightarrow{e_1}311|3311|221|3\]
is computable, but not perhaps obvious. However, represented as perforated tableaux, the result of successive applications of $e_1$ is easy to see:
\[^*T=\begin{ytableau} \ & &&2&2&\\ &1&1&3&3&3\\1 &2&2&4&& \end{ytableau}\xrightarrow{e_1}\begin{ytableau} \ & &&2&2&3\\ &1&1&3&3&\\1 &2&2&4&& \end{ytableau}\xrightarrow{e_1}\begin{ytableau} \ & &1&2&2&3\\ &1&&3&3&\\1 &2&2&4&& \end{ytableau}\xrightarrow{e_1}\begin{ytableau} \ &1 &1&2&2&3\\ &&&3&3&\\1 &2&2&4&& \end{ytableau},\]
and the action of $e_1$ becomes even more evident when restricting our attention to $^*T[1,2]$, where applying $e_1$ just moves left-justified content under a blank up:

\[ ^*T[1,2]=\begin{ytableau} \ &&2&2&\\ 1&1&3&3&3\end{ytableau}\xrightarrow{e_1}\begin{ytableau} \ &&2&2&3\\ 1&1&3&3&\end{ytableau} \xrightarrow{e_1}\begin{ytableau} \ &1&2&2&3\\ 1&&3&3&\end{ytableau}\xrightarrow{e_1}\begin{ytableau} 1 &1&2&2&3\\ &&3&3&\end{ytableau}.\]

Below are, we claim, two (isomorphic) crystal graphs of ptableaux. The crystal on the left is the image under $\Pf$ of a crystal in $[3]^{\otimes 4}$ with highest weight $1112$. This crystal can be realized explicitly as a crystal of SSYT under a row-reading map. However, the crystal on the right corresponds to highest weight $1211 \in [3]^{\otimes 4}$, which does not have a corresponding SSYT model.

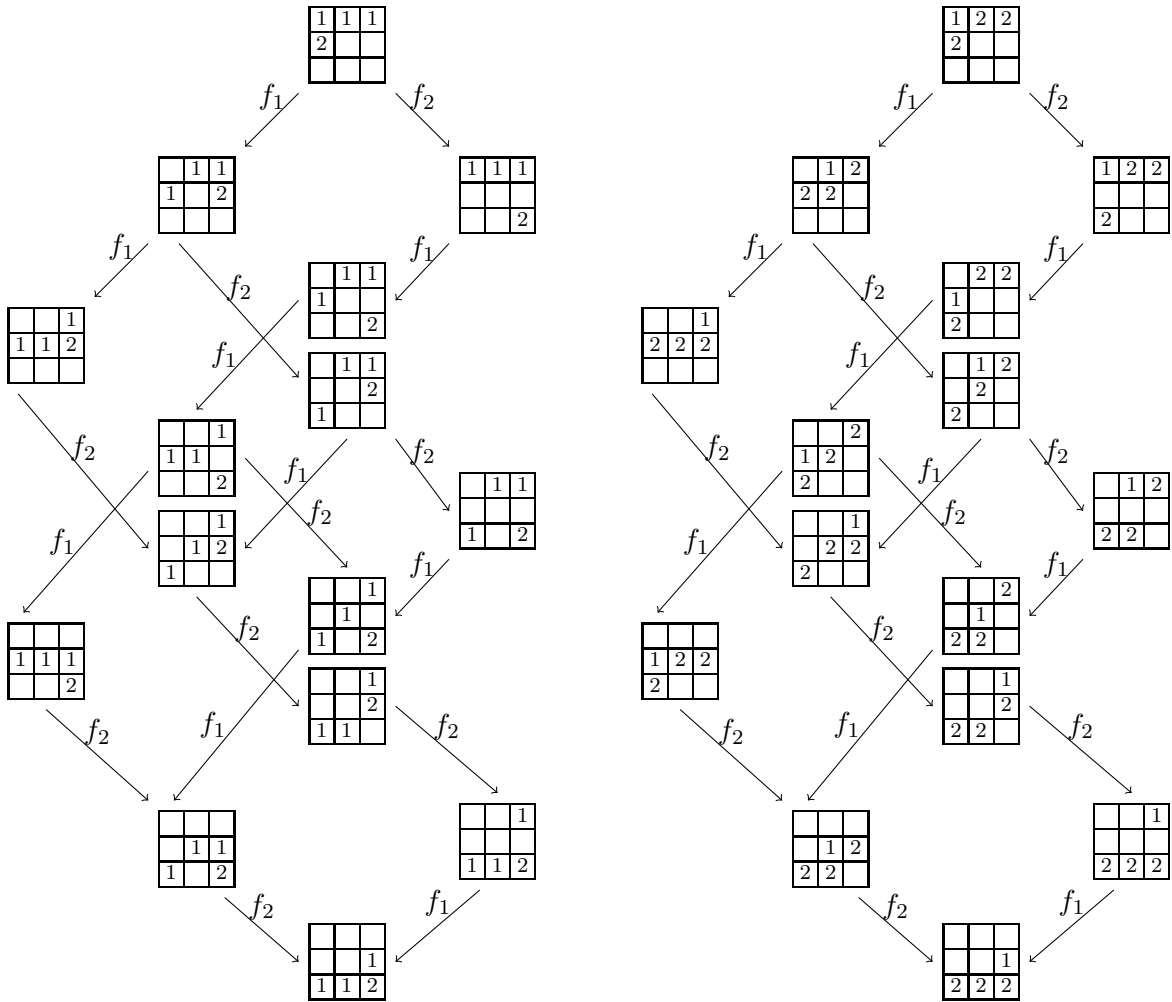
\[ \ytableausetup{smalltableaux}
\begin{tikzpicture}
\node (A) { \ytableaushort{111,2,}* {3,3,3}};
\node [below of=A, left of=A,node distance = 2cm](B) {\ytableaushort{\none11,1\none2,}* {3,3,3}};
\node (C) [right of = A,below of=A,node distance = 2cm] {\ytableaushort{111,\none,\none\none2}* {3,3,3}};
\node (D) [left of=B, below of = B,node distance = 2cm] {\ytableaushort{\none\none1,112,}* {3,3,3}};
\node (E) [below of=A,node distance = 3.4cm] {\ytableaushort{\none11,1,\none\none2}* {3,3,3} };
\node (F) [below of=E,node distance =1.2cm] {\ytableaushort{\none11,\none\none2,1}* {3,3,3} };
\node (G) [below of =B,node distance= 3.5cm]  {\ytableaushort{\none\none1,11,\none\none2}* {3,3,3}};
\node (H) [below of = G,node distance = 1.2 cm] {\ytableaushort{\none\none1,\none12,1}* {3,3,3} };
\node (I) [below of =C,node distance = 4.2cm] {\ytableaushort{\none11,\none\none\none,1\none2}* {3,3,3}};
\node (J) [below of=D, node distance = 4.2cm] {\ytableaushort{\none\none\none,111,\none\none2}* {3,3,3}};
\node (K) [below of = F,node distance = 3cm] {\ytableaushort{\none\none1,\none1,1\none2}* {3,3,3}  };
\node (L) [below of = K,node distance = 1.2 cm] {\ytableaushort{\none\none1,\none\none2,11}* {3,3,3}};
\node (M) [below of = H,node distance = 4 cm] { \ytableaushort{\none\none\none,\none11,1\none2}* {3,3,3}};
\node (N) [below of = I,node distance = 4.4 cm]{\ytableaushort{\none\none1,\none\none\none,112}* {3,3,3}};
\node (O) [below of = L,node distance = 3.4cm] {\ytableaushort{\none\none\none,\none\none1,112}* {3,3,3}};
\draw[->] (A.225) to node [above] {$f_{1}$} (B.45);
\draw[->] (A.315) to node [above] {$f_{2}$} (C.135);
\draw[->] (B.225) to node [above] {$f_{1}$} (D.45);
\draw[->] (B.250) to node [above] {$f_{2}$} (F.165);
\draw[->] (C.225) to node [above] {$f_{1}$} (E.0);
\draw[->] (E.180) to node [left] {$f_{1}$} (G.90);
\draw[->] (D.240) to node [above] {$f_{2}$} (H.180);
\draw[->] (G.195) to node [left] {$f_{1}$} (J.115);
\draw[->] (F.270) to node [above] {$f_{1}$} (H.0);
\draw[->] (F.315) to node [above] {$f_{2}$} (I.180);
\draw[->] (H.270) to node [above] {$f_{2}$} (L.180);
\draw[->] (I.225) to node [above] {$f_{1}$} (K.0);
\draw[->] (J.270) to node [above] {$f_{2}$} (M.135);
\draw[->] (G.0) to node [right] {$f_{2}$} (K.90);
\draw[->] (M.300) to node [above] {$f_{2}$} (O.180);
\draw[->] (K.215) to node [left] {$f_{1}$} (M.115);
\draw[->] (N.250) to node [above] {$f_{1}$} (O.0);
\draw[->] (L.0) to node [above] {$f_{2}$} (N.90);
\end{tikzpicture} \qquad \quad
\begin{tikzpicture}
\node (A) { \ytableaushort{122,2,}* {3,3,3}};
\node [below of=A, left of=A,node distance = 2cm](B) {\ytableaushort{\none 1 2,22,}* {3,3,3}};
\node (C) [right of = A,below of=A,node distance = 2cm] {\ytableaushort{122,\none,2}* {3,3,3}};
\node (D) [left of=B, below of = B,node distance = 2cm] {\ytableaushort{\none\none1,222,}* {3,3,3}};
\node (E) [below of=A,node distance = 3.4cm] {\ytableaushort{\none 22,1 \none \none,2}* {3,3,3} };
\node (F) [below of=E,node distance =1.2cm] {\ytableaushort{ \none1 2,\none 2 \none,2}* {3,3,3} };
\node (G) [below of =B,node distance= 3.5cm]  {\ytableaushort{\none \none 2 , 12,2 }* {3,3,3}};
\node (H) [below of = G,node distance = 1.2 cm] {\ytableaushort{\none\none1,\none 2 2,2}* {3,3,3} };
\node (I) [below of =C,node distance = 4.2cm] {\ytableaushort{ \none1 2,\none ,22}* {3,3,3}};
\node (J) [below of=D, node distance = 4.2cm] {\ytableaushort{\none,122,2}* {3,3,3}};
\node (K) [below of = F,node distance = 3cm] {\ytableaushort{\none \none 2,\none1,22}* {3,3,3}  };
\node (L) [below of = K,node distance = 1.2 cm] {\ytableaushort{\none \none1, \none \none 2,22}* {3,3,3}};
\node (M) [below of = H,node distance = 4 cm] { \ytableaushort{\none, \none 1 2,22}* {3,3,3}};
\node (N) [below of = I,node distance = 4.4 cm]{\ytableaushort{\none \none1,\none,222}* {3,3,3}};
\node (O) [below of = L,node distance = 3.4cm] {\ytableaushort{\none,\none\none1,222}* {3,3,3}};
\draw[->] (A.225) to node [above] {$f_{1}$} (B.45);
\draw[->] (A.315) to node [above] {$f_{2}$} (C.135);
\draw[->] (B.225) to node [above] {$f_{1}$} (D.45);
\draw[->] (B.250) to node [above] {$f_{2}$} (F.165);
\draw[->] (C.225) to node [above] {$f_{1}$} (E.0);
\draw[->] (E.180) to node [left] {$f_{1}$} (G.90);
\draw[->] (D.240) to node [above] {$f_{2}$} (H.180);
\draw[->] (G.195) to node [left] {$f_{1}$} (J.115);
\draw[->] (F.270) to node [above] {$f_{1}$} (H.0);
\draw[->] (F.315) to node [above] {$f_{2}$} (I.180);
\draw[->] (H.270) to node [above] {$f_{2}$} (L.180);
\draw[->] (I.225) to node [above] {$f_{1}$} (K.0);
\draw[->] (J.270) to node [above] {$f_{2}$} (M.135);
\draw[->] (G.0) to node [right] {$f_{2}$} (K.90);
\draw[->] (M.300) to node [above] {$f_{2}$} (O.180);
\draw[->] (K.215) to node [left] {$f_{1}$} (M.115);
\draw[->] (N.250) to node [above] {$f_{1}$} (O.0);
\draw[->] (L.0) to node [above] {$f_{2}$} (N.90);
\end{tikzpicture}
 \]

\section{First Results}
The following theorem will imply that the bijections mentioned in the introduction, and proved in Lemma~\ref{Perf order} along with the definitions for crystal operators on ptableaux of the previous section, induce a $GL_n$ crystal graph structure on ptableaux. While one could simply \emph{define} crystal operators on ptableaux by transfer of structure~\cite{BumpSchilling,Shimozono}, we would still need to verify that this transfer is actually realized by the simpler {jeu de taquin} style maps constructed above, saving ourselves little work, and so we present the theorem below.
\begin{thm}
Let $e_{i}: P([n]^{\otimes k}) \rightarrow P([n]^{\otimes k}) \cup \{ \NULL \}$ be the standard crystal operator and let $\omega_{\cal P} \in P([n]^{\otimes k}) $, with some parsing ${\cal P}$. Then,
\[ e_{i} \left( \Pf( \omega_{\cal P}) \right) = \Pf\left(e_{i}(\omega)_{\cal P} \right). \] \label{ei map}
That is, the standard $e_{i}$ crystal operator on words in $[n]^{\otimes k}$ commutes with the $e_{i}$ operator on ptableaux (Definition~\ref{e_i def}).
\end{thm}

The above equality also implies that if $\omega_{\cal P}$ is some parsing of a word $\omega \in [n]^{\otimes k} $, then ${\cal P}$ remains a valid parsing of $e_{i}(\omega)$, so that we can write $e_{i}(\omega_{\cal P}) = e_{i}(\omega)_{\cal P}$.

The effect of $e_{i}$ on some word $\omega_{\cal P}$ (if the output is not $\NULL$), is to change a letter $(i+1)$ in $\omega_{\cal P}$ into an $i$ in $e_{i}(\omega_{\cal P})$. Thus, the difference between $\Pf(\omega_{\cal P})$ and $\Pf(e_{i}(\omega_{\cal P}))$ is that exactly one entry in row $i+1$ of $\Pf(\omega_{\cal P})$ is moved to row $i$ in $\Pf(e_{i}(\omega_{\cal P}))$. Theorem~\ref{ei map} says that our definition of $e_{i}$ on ptableaux chooses the correct entry to move, by swapping it with a blank lying directly above it.

Before proving the theorem, we give an example. Given some parsed word $\omega_{\cal P}$, let $\omega_{\cal P}[i,i+1]$ denote the sequence of $i$'s and $(i+1)$'s appearing in $\omega_{\cal P}$, in order, maintaining the parsing. So, suppose, for some word, we had:

\[ \omega_{\cal P}[3,4]= 44433|444333|4433|44443|433. \]

Recall the definition of $c_{e}^{i}(\omega,j)$ Theorem~\ref{crystal def}. We write the running totals of $c_{e}^3(\omega_{\cal P},j)$ (from left to right) under the entries in $\omega_{\cal P}$, putting the occurrence of any \emph{new} maximum in bold:
\begin{align*}  \omega_{\cal P}[3,4]= &
\mathbf{444}33|44\mathbf{4}333|4433|444\mathbf{4}3|4\mathbf{4}3\\
&\mathbf{123}32|23\mathbf{4}432|2332|234\mathbf{5}5|5\mathbf{6}6
 \end{align*}

Let $T=\Pf(\omega_{\cal P})$.  $T$ is constructed, under $\Pf$, by reading the letters of $\omega_{\cal P}$, left to right. We say a box $\framebox{$a$}$ in $T$ is ``read into" the ptableau ``before" some other box $\framebox{$b$}$ whenever $\Pf^{-1}(\,\framebox{$a$}\,)$ lies to the left of $\Pf^{-1}(\,\framebox{$b$}\,)$ in $\omega_{\cal P}$. We label the entries of $^*T[i,i+1]$ with subscripts, indicating the order (from left to right) that they are read into the ptableau:
\ytableausetup{nosmalltableaux}
\[ ^*T[3,4] = \begin{ytableau}\  &&&\scalebox{0.6}[0.7]{$1_4$} & \scalebox{0.6}[0.7]{$1_5$} &&
 \scalebox{0.6}[0.7]{$2_9$} &\scalebox{0.6}[0.7]{$2_{10}$}&\scalebox{0.6}[0.7]{$2_{11}$} &\scalebox{0.6}[0.7]{$3_{14}$} & \scalebox{0.6}[0.7]{$3_{15}$} & & \scalebox{0.6}[0.7]{$4_{20}$}&&\scalebox{0.6}[0.7]{$5_{23}$} \\
\scalebox{0.6}[0.7]{$\mathbf{1_1}$} & \scalebox{0.6}[0.7]{$\mathbf{1_2}$}& \scalebox{0.6}[0.7]{$\mathbf{1_3}$} & \scalebox{0.6}[0.7]{$2_6$} & \scalebox{0.6}[0.7]{$2_7$} & \scalebox{0.6}[0.7]{$\mathbf{2_8}$}&\scalebox{0.6}[0.7]{$3_{12}$}& \scalebox{0.6}[0.7]{$3_{13}$} & \scalebox{0.6}[0.7]{$4_{16}$}& \scalebox{0.6}[0.7]{$4_{17}$}&\scalebox{0.6}[0.7]{$4_{18}$}& \scalebox{0.6}[0.7]{$\mathbf{4_{19}}$} & \scalebox{0.6}[0.7]{$5_{21}$}&\scalebox{0.6}[0.7]{$\mathbf{5_{22}} $} &\end{ytableau} \]
We see that each new maximum for $c_{e}^3(\omega_{\cal P},j)$ corresponds, under $\Pf$, to  an entry in $^*T[3,4]$ in row $4$ with no content above it. Thus, $e_{3}(\omega_{\cal P})$ changes the rightmost $4$ to a $3$, while $e_{2}(^*T[3,4])$ (and hence, $e_{2}(T)$ as well) moves the corresponding entry ($5_{22}$, in the penultimate column) from row $4$ to row $3$. Note also that the content to be moved is immediately apparent from the form $^*T[3,4]$, without the computation needed in $\omega_{\cal P}[3,4]$.

\bigskip

\ytableausetup{smalltableaux}

\begin{proof}
 Let
$\omega = \omega_{1}  \omega_{2} \cdots \omega_{k-1} \omega_k \in [n]^{\otimes k}$
and choose some parsing $\omega_{\cal P} = h_{1}|h_{2}| \cdots |h_s. $
We show $e_{i}( \Pf(\omega_{\cal P})) = \Pf(e_{i}(\omega)_{\cal P})$, and that ${\cal P}$ is still a valid parsing of $e_{i}(\omega)$. Let $T = T_{\omega_{\cal P}} = \Pf(\omega_{\cal P})$.

It is clear that
$^*T[i,i+1] = \Pf(\omega_{\cal P}[i,i+1]).$
Further, since $^*T[i,i+1]$ has only two rows, the bottom row (row $i+1$), is, from left to right, a block of content followed by a block of blanks. In row $i$, content is as far left as possible, so that, if any box  $\framebox{$c'$}$ has a blank to its left, it is to maintain column-strictness.  Thus, in the diagram below, $c\le c'$ and $c''\le c'''$ with blanks above the entries $\framebox{$c$}$ and $\framebox{$c''$}$:

\[^*T[i,i+1] = \begin{ytableau} \none & \none & \none & c' & \cdots & \cdots & \none & \none & c''' & \cdots & \cdots\\ \cdots & \cdots & c & a & \cdots & \cdots & \cdots & c''& b& \cdots \end{ytableau} \]

The bijection induced from $\Pf$ taking letters in $\omega_{\cal P}$ to boxes in $T_{\omega_{\cal P}}$ can be restricted to $\Pf:\omega_{\cal P}[i,i+1] \rightarrow T[i,i+1]$.
Suppose $\framebox{c}$ in row $i+1$ of $^*T[i,i+1]$ is an \emph{uncovered entry}, meaning it has a blank directly above it. Let $x=\Pf^{-1}(\,\framebox{c}\,)$. Since $^*T[i,i+1]$ is left-justified, and column-strict, the letters to the left of $x$ in $\omega_{\cal P}[i,i+1]$ map bijectively, under \Pf, to the entries left of $\framebox{$c$}$ in rows $i$ or $i+1$ of $^*T[i,i+1]$. Then (again recalling the definition of $c_{e}^{i}(\omega,j)$ from Theorem~\ref{crystal def}):
\begin{align*} 0 &< ( \hbox{number of entries in row $(i+1)$  weakly left of $\framebox{c}$\,} ) -\
\!\!(\hbox{number of entries in row $i$ left of $\framebox{c}$\,} ) \\
&= \ \!\! ( \hbox{number of $(i+1)$'s weakly left of $x$} ) -\
\!\!( \hbox{number of $i$'s left of $x$} ) \\
&= c_{e}^{i}(\omega_{\cal P},j),\end{align*}
where we assume $x$ appears in column $j$ of $\omega_{\cal P}$.

Thus, uncovered entries in row $i+1$ of $^*T[i,i+1]$ correspond to positive values of $c_{e}^{i}(\omega_{\cal P},j)$, with the largest such value corresponding to the \emph{rightmost} uncovered entry.

We claim only the \emph{first} occurrence of a new maximum value of $c_{e}^{i}(\omega_{\cal P},j)$ (as we compute values, left to right, across the word $\omega_{\cal P}$) will correspond to an uncovered entry of $^*T[i,i+1]$ in row $(i+1)$. Indeed, if an entry $\framebox{b}$ in row $i+1$ is \emph{covered}, then the entries directly above or to its left will have been read in prior to it, (although some entries that were read in before
$\framebox{b}$ may lie in row $i$ to the right of $\framebox{b}$). Thus, the value of $c_{e}^{i}(\omega_{\cal P},j')$ corresponding to $\framebox{b}$ either cannot exceed that of an uncovered entry to its immediate left, or is negative,(if there is no such uncovered entry).  Therefore the first occurrence of the maximal value of $c_{e}^{i}(\omega_{\cal P},j)$ cannot occur at $\framebox{b}$.

Thus we see that the first appearance of the maximum value of $c_{e}^{i}(\omega_{\cal P},j)$ can only appear in an entry corresponding to the \emph{rightmost} uncovered entry of $^*T[i,i+1]$ in row $(i+1)$.
This rightmost uncovered entry will be moved into row $i$ in the computation of $e_{i}(T)$ precisely when the corresponding entry in $\omega_{\cal P}$ (under \Pf) changes an $i+1$ to an $i$ in computing $e_{i}(\omega_{\cal P})$, proving the result.
\end{proof}

\begin{df} Define a map $\Rot : \PTab_{(\ell,n)} \rightarrow \PTab_{(\ell,n)}$ (recall the definition of $\PTab_{(\ell,n)}$ from Definition~\ref{ptab def}) with $\Rot(T)$ the ptableau obtained from $T$ by first rotating $T$ by $180^\circ$, and then replacing each entry $t$ by $\ell-t+1$.
\end{df}

For example, if
\ytableausetup{smalltableaux}
 $T =  \begin{ytableau}
\ &&&&&&&&4&4&5\\
&&&&&1&1&2&&5&6\\
&&&1&1&&2&4&5&6&7\\
1&1&2&3&3&3&4&6&6&&8 \end{ytableau} \in PTab_{(8,4)}, $

\vskip 0.1in

\noindent then
$ \Rot(T) = \begin{ytableau}
\ 1 & & 3 &3&5&6&6&6&7&8&8 \\
2&3&4&5&7&&8&8&&&\\
3&4&&7&8& 8&&&&&\\
4&5&5&&&&&&&& \end{ytableau}\, . $

Similarly,
\begin{df} Define a map $\Rot : [n]^{\otimes k} \rightarrow  [n]^{\otimes k}$ by first defining, for each $\ell \in [n]$:
\[ \overline{\ell} = n- \ell +1, \]
and then defining
\[ \Rot(\omega) = \Rot(\omega_{1}\omega_{2} \dots \omega_{k}) = \overline{\omega_{k}}\,\overline{\omega_{k-1}} \dots \overline{\omega_{1}}. \]
\end{df}

We state the following easy facts, without proof:
\begin{prop}
For any $T \in \PTab_{n}$, we have $^{*}\Rot(T)=\Rot(T^{*})$ and $\Rot(^{*}T)=\Rot(T)^{*}$. Also, the map $\Rot$ is a involution, and finally
\[ \Pf(\Rot(\omega)) = \Rot(\Pf(\omega)). \]
\end{prop}

Analogously to Definition \ref {e_i def}, we define the ptableau crystal operator $f_i$ ($1\le i < n$).  Note that a result for right-justified ptableaux, analogous to Lemma 4.6, ensures that the following is well-defined:

\begin{df}  \label {f_i def} Suppose $T \in \PTab_{n}$ and that $\framebox{c}$ is an entry in row $i$ of $T$ such that in $T[i,i+1]^*$ (the right-justified form of $T[i,i+1]$) $\framebox{c}$ is the left-most entry with a blank below it. In this case, we define $f_i (T)$ to be the ptableau obtained by swapping $\framebox{c}$ with the blank below it (in row $i+1$) of $T^*$. If there is \emph{no} such entry $\framebox{c}$, then $f_i (T) = NULL$.
\medskip

\noindent We define $\phi_{i}(T)$ to be the number of blanks in row $i+1$ of $T[i,i+1]^*$.
\end{df}

\begin{thm} Given the map $f_{i} : \PTab_{n} \rightarrow \PTab_{n}$ above, we have
\[ \Pf(f_{i}(\omega_{\cal P})) = f_{i}(\Pf(\omega_{\cal P})). \]
\end{thm}
\begin{proof} We note (from the definitions of the maps $e_{i}$ and $f_{i}$) that $\Rot(f_{i}(\omega_{\cal P})) = e_{n-i}(\Rot(\omega_{\cal P}))$ and
$\Rot(e_{n-i} (T)) = f_{i}(\Rot(T))$, from which we have
\begin{align*} \Pf(f_{i}(\omega_{\cal P})) & = \Pf(\Rot(\Rot(f_{i}(\omega_{\cal P})))) \\
& = \Pf(\Rot(e_{n-i}(\Rot(\omega_{\cal P})))) \\
& = \Rot(\Pf (e_{n-i}(\Rot(\omega_{\cal P}))) \\
& = \Rot(e_{n-i}(\Pf(\Rot(\omega_{\cal P})))) \qquad \hbox{\emph{by Theorem~\ref{ei map}}} \\
& = f_{i}(\Rot(\Pf(\Rot(\omega_{\cal P})))) \\
& = f_{i}(\Pf(\Rot(\Rot(\omega_{\cal P})))) \\
& = f_{i} (\Pf(\omega_{\cal P})). \end{align*}\end{proof}

The following Corollary highlights the convenience of computing crystal operators in $\PTab_n$.  Unlike the complicated calculations for $\phi_i$ and $\epsilon_i$ for words in $[n]^{\otimes k}$, we can calculate the corresponding values for ptableaux with a simple count in $T[i,i+1]$, without needing left- or right-justification:

\begin{cor} Suppose $T \in PTab_n$. Then we have
\[ \phi_i(T) = \hbox{the number of blanks in row $i+1$ of $T[i,i+1]$} , \]
and
\[ \epsilon_i(T) =   \hbox{the number of blanks in row $i$ of $T[i,i+1]$} . \]\label{ep phi def}
\end{cor}
The reader should compare the simplicity of the definitions for $\phi_i(T), \epsilon_i(T)$ with the corresponding definitions of $\phi_i(\omega)$ and $\epsilon_i(\omega)$ in Proposition~\ref{crystal def}.
\medskip

\begin{cor} Let $\omega_{\cal P}$ be highest weight, let $T_{\omega_{\cal P}} = \Pf(\omega_{\cal P})$, and $\nu = \wt(\omega_{\cal P}) = \wt(T_{\omega_{\cal P}})$.  Let ${\mathbb B}^{\nu}(T_{\omega_{\cal P}} )$ be the closure of the images of (highest weight) $T_{\omega_{\cal P}}$ under the ptableaux operators $e_{i}$ and $f_{i}$, for $1 \leq i < n$. Then the restriction $\Pf: {\mathbb B}^{\nu}(\omega_{\cal P}) \rightarrow {\mathbb B}^{\nu}(T_{\omega_{\cal P}} )$ is an isomorphism of irreducible $GL_n$ crystals. (Recall the definition of ${\mathbb B}^{\nu}(\omega_{\cal P})$ from Definition~\ref{word irred}.)
\end{cor}
\begin{proof} The set ${\mathbb B}^{\nu}(\omega_{\cal P})$ is a connected, irreducible subcrystal of $[n]^{\otimes k}$. Since $\Pf$ is a bijection from ${\mathbb B}^{\nu}(\omega_{\cal P})$ that commutes with the crystal operators, and $\wt(T_{\omega_{\cal P}}) = \wt(\omega_{\cal P})$ as noted above, they are isomorphic crystals (See~\cite{BumpSchilling}, p. 54).
\end{proof}
\medskip

\section{Highest Weights}
\begin{df} We say a ptableau $T \in \PTab$ is \emph{partition-shaped} if no blank box in $^{*}T$ has content to the right or below the blank (so all content is in the upper left corner of $^{*}T$, with no ``holes").
 We say a ptableau is \emph{anti-partition-shaped} if in $T^{*}$, no blank has content to the left or above it or, equivalently, if $\Rot(T)$ is partition-shaped.
\end{df}

For example, the ptableau
\[ T = \begin{ytableau}
\ &&&&1&1&4 \\
&&1&1&&2&5\\
1&2&3&4&4&5& 6\end{ytableau} \]
is not partition-shaped since
\[ ^{*}T =  \begin{ytableau}
\ &&&1&1&4 &\\
&1&1&2&&&5\\
1&2&3&4&4&5&6 \end{ytableau}, \]
and the leftmost blanks in the first two rows have content to their right. However,
the ptableau
\[ T' =  \begin{ytableau}
 1&1&1&1&1&4&4 \\
& 2 && 4&4&6&\\
&&3&&5&& \end{ytableau}, \]
is partition-shaped since
\[ ^{*}T' = \begin{ytableau}
 1&1&1&1&1&4&4 \\
 2 & 4&4&6&&&\\
3&5&&&&& \end{ytableau} . \]

\begin{thm} A ptableau $T \in \PTab_{n}$ is partition-shaped if and only if it is highest weight, meaning $e_{i}(T) = \NULL$ for all $i$, $1 \leq i <n$.

Similarly, a ptableau $T$ is anti-partition-shaped if and only if $T$ is \emph{lowest weight}, meaning $f_{i}(T) = \NULL$ for all $i$, $1 \leq i <n$.
\label{highest weight}
\end{thm}
Thus, in the examples above, the ptableau $T'$ is highest weight, while $T$ is not.
\medskip

\begin{proof} Clearly, if $^{*}T$ is partition-shaped, then it is highest weight, since $e_{i}$ on $^{*}T$, (if not $\NULL$), swaps a blank in row $i$ with content directly below it in row $(i+1)$.  When $^{*}T$ is partition-shaped, there are no such blanks. Equivalently, by Corollary~\ref{ep phi def} all $\epsilon_i(T)=0$.

For the converse, we note that by Theorem~\ref{ei map} and Definition~\ref{e_i def}, if $e_i (T) = NULL$ for all $i$, each $^{*}T[i,i+1]$ must be partition-shaped; otherwise there would be some blank over content, contradicting $e_i (T) = NULL$. Hence $T$ itself is partition shaped.

The proof of the second statement regarding lowest weights is proved analogously (indeed, isomorphically under the map $\Rot$). \end{proof}

It is well-known that a word $\omega \in [n]^{\otimes k}$ is a highest weight if and only if it is a \emph{Yamanouchi} word: if $\omega = \omega_{1} \omega_{2} \cdots \omega_{k-1} \omega_{k}$, then for each $i \in [n-1]$ and each $ s$, $1 \leq s \leq k$, the segment $\omega_{1} \omega_{2} \cdots \omega_{s-1} \omega_{s}$ contains at least as many $i$'s as $(i+1)$'s. The characterization of highest weight ptableau given above is, by contrast, generally easier to use. For example, to see if the minimally parsed word
\[ \omega= 111|22111|2221|33211|332|33 \]
is highest weight, we could check that all $e_{i}$ are $\NULL$ on $\omega$, or we could compute all the counts of $i$'s over $(i+1)$'s across the length of the word to verify that the word is Yamanouchi, or, much more simply, we can make
the associated ptableau:
\ytableausetup{smalltableaux}
 \ytableausetup{centertableaux}
\[ \Pf(\omega) = \begin{ytableau} 1 & 1 & 1& 2 & 2&2& 3 & 4 &4 \\
2 & 2& 3 & 3& 3& 4 & 5& & \\
4 & 4& 5 & 5 & 6 & 6& && \end{ytableau}, \]
and observe immediately that it is highest weight (unperforated). Indeed, viewing the $e_{i}$ operators as moving blanks down to lower rows in a ptableau, and the $f_{i}$ as moving blanks up, it is quite natural that highest weight ptableaux have no interior blanks; the blanks have been moved down as far as possible.

\section{Tensor Products and the Littlewood-Richardson Rule}

\subsection{Tensor Products in $[n]^{\otimes k}$}
Tensor products of crystals in $[n]^{\otimes k}$ are constructed by concatenation of words. If ${\cal C}_{1}$ and ${\cal C}_{2}$ are two crystal graphs (from, say, $[n]^{\otimes k}$ and $[n]^{\otimes s}$, respectively), then the nodes of the crystal graph ${\cal C}_{1} \otimes {\cal C}_{2} \subseteq [n]^{\otimes (k + s)}$ are the concatenated words of ${\cal C}_{1}$ and ${\cal C}_{2}$, with the crystal operator acting on these words.   The set of nodes in ${\cal C}_{1} \otimes {\cal C}_{2}$ is \emph{closed} under the crystal operators, and hence decomposes as a disjoint sum of irreducible crystals (see \cite{BumpSchilling}, p.\ 18).

Let ${\mathbb B}^{\nu}(\omega)$ and ${\mathbb B}^{\mu}(\omega')$ be two irreducible crystals in $[n]^{\otimes k}$ with highest weight elments $\omega$ and $\omega'$, respectively, with weights $wt(\omega) = \nu$ and $wt(\omega') = \mu$. A classical problem is to determine the number of irreducible crystal graphs of a given highest weight $\lambda$ in the tensor product ${\mathbb B}^{\nu}(\omega) \otimes {\mathbb B}^{\mu}(\omega')$. We denote this multiplicity by the \emph{Littlewood-Richardson coefficient}, $c_{\mu \nu}^{\lambda}$. In order to determine $c_{\mu \nu}^{\lambda}$ \emph{combinatorially}, it is not enough to use representatives of isomorphism classes, such as SSYT.

The issue is the SSYT crystal model is not closed under tensor products. As noted above, if $\nu$ is a partition of $n$ rows and $k$ parts, the map $RR: SSYT(\nu) \rightarrow [n]^{\otimes k}$ picks out only \emph{one} representative irreducible crystal of highest weight $\nu$, with all other graphs with the same highest weight identified only up to isomorphism. So, for example, the pair of SSYT:
\[ T = \begin{ytableau} 1 & 1 & 3 \\ 2 & 3 \end{ytableau} \quad \hbox{and} \quad T' = \begin{ytableau} 1 & 2 & 3 \\ 2 & 3 \end{ytableau} \]
correspond, under the row-reading map $RR$ to the words:
\[ \begin{array}{ccccc}

T = \begin{ytableau} 1 & 1 & 3 \\ 2 & 3 \end{ytableau} & & T' = \begin{ytableau} 1 & 2 & 3 \\ 2 & 3 \end{ytableau}& & \\
 & & & & \\
\Downarrow & & \Downarrow && \\
 & & & & \\
RR(T) = 31132  & & RR(T') = 32132 & & \end{array}. \]

We can certainly form the tensor product of the words:
\[ RR(T) \otimes RR(T') =  31132 \otimes 32132 = 31132|32132, \]
(the parsing is there for clarity, and does not affect the crystal structure). However, this word has weakly decreasing factors (from left to right) of size $3,2,3,2$, and hence is clearly \emph{not} the row-reading word of any SSYT (for which the lengths of factors must be weakly decreasing). Of course, there are tableau methods to replace $31132|32132$ with a the row-reading word of some SSYT, but this requires computing an isomorphic \emph{proxy} for this tensor product.

To compute \LR\ coefficients \emph{combinatorially}, one must employ a different combinatorial model.
Classically, this is achieved by enumerating \emph{Littlewood-Richardson fillings} of the skew shape $\lambda / \mu$ with content $\nu$. Hence, in the setting of $[n]^{\otimes k}$, we use one class of objects (semistandard Young tableau) to build representative models of crystals, but to analyze tensor products, we need to build a different model (fillings of skew shapes). Using ptableaux, we can accomplish \emph{both} tasks simultaneously.

\subsection{Tensor Products of Ptableaux}

\begin{df}We adopt the following convention. Given two words $\omega, \omega' \in [n]^{\otimes k}$, with parsings ${\cal P}$ and ${\cal P}'$, resulting in $\omega_{\cal P}$ and $\omega_{{\cal P'}}'$, we will assign to $\omega_{\cal P} \otimes w_{{\cal P'}}'$ the parsing
\[ \omega_{\cal P} \otimes \omega_{{\cal P'}}' = \omega_{\cal P}| \omega_{{\cal P'}}'. \]
\end{df}

The tensor product of ptableaux is defined quite easily. Let $T \in \PTab_{(\ell,n)}$ and $U \in \PTab_{(m,n)}$. We form the \emph{concatenation} $T \otimes U \in\PTab_{((\ell+m),n)}$ by:
\begin{enumerate}
\item Replace each box $\framebox{$a$}$ in $U$ by $\framebox{$a + \ell$}$. Call this ptableau $U^{+\ell}$.
\item Construct $T \otimes U$ by appending the rows of $U^{+\ell}$ to the \emph{right} of the rows of $^{*}T$, from the bottom to the top, at each stage left-justifying appended rows. For example, if
\end{enumerate}

$T = \begin{ytableau}
\ & 1 & 1 & 2 & & 3 \\ 1 & 2 & 3 & 3 &3 &   \\ 2&3 & 4 & && \end{ytableau}, U = \begin{ytableau} 1 & 1 & 1 & 2 & 4 & 4  \\ 2 & 2 & & 3 & & \\ 3 & 3 & 3 &  & &  \end{ytableau}\in \PTab{(4,3)}$,
then
$  U^{+4} = \begin{ytableau} 5 & 5 & 5 & 6 & 8 & 8  \\ 6 & 6 & & 7 & & \\ 7 & 7 & 7 &  & &  \end{ytableau}\in \PTab{(8,3)}$,
and so
\[ T \otimes U = \begin{ytableau}
\ &1 & 1 & 2 & & 3 & 5 & 5 & 5 & 6 & 8 & 8 \\
1 & 2 & 3 & 3 & 3 & 6 & 6 & 7 &&& &\\ 2 & 3 & 4 & 7 & 7 & 7 & &&&&& \end{ytableau}\, . \]

Note that it is possible that a blank appearing in $U$ will not persist in $T \otimes U$.

The following is easily proved:
\begin{lem} Concatenation of words in $[n]^{\otimes k}$ commutes with $\Pf$ in $\PTab_{n}$. That is,
\[ \Pf (\omega_{\cal P} \otimes \omega'_{\cal P'}) = \Pf (\omega_{\cal P}) \otimes \Pf(\omega'_{\cal P'}). \] \label{tensor parse}
\end{lem}
As noted earlier, a word $ \omega \in [n]^{\otimes k}$ is of highest weight if and only if it is a \emph{Yamanouchi} word, that is,  if $\omega = \omega_{1} \omega_{2} \cdots \omega_{k-1} \omega_{k}$, then for each $i \in [n-1]$ and each $ \ell$, $1 \leq s \leq k$, the segment $\omega_{1} \omega_{2} \cdots \omega_{s -1} \omega_{s}$ contains at least as many $i$'s as $(i+1)$'s. Thus, a concatenation of words $\omega \otimes \omega'$ is highest weight precisely when it is Yamanouchi. In the ptableau setting, this is essentially a corollary to Theorem~\ref{highest weight}:

\begin{thm} Let $T,T' \in \PTab_{n}$. Then $T \otimes T'$ is highest weight if and only if $T \otimes T'$ is partition-shaped. In particular, this implies that $T$ is highest weight (and hence itself partition-shaped).
\end{thm}
\begin{proof} By Theorem~\ref{highest weight}, $T \otimes T'$ is highest weight if and only if $T \otimes T'$ is partition-shaped. Since $T$ forms the interior, upper-left corner of $T \otimes T'$, $T \otimes T'$ will be highest weight if and only if $T$ already has no interior blanks (in its left-justified form), in other words, $T$ itself is highest weight.
\end{proof}

For example, suppose
$ T =\begin{ytableau}
1 & 1 & 1 & 2 & 2& 3& 4 & 4 \\ 1 & 2 & 3 & 3 &3 &&&   \\ 2&3 & 4 & &&&& \end{ytableau}$
and
$T' = \begin{ytableau}
\ & 1 & 2& 2 & 2& 3 \\ 1 & 2 &  &  3&3 &   \\ 2&3 & 3 & && \end{ytableau} $.

Then
\[ T \otimes T' =  \begin{ytableau}
1 & 1 & 1 & 2 & 2& 3& 4 & 4 & 5 & 6& 6& 6& 7 \\ 1 & 2 & 3 & 3 &3 &5 &6 &    7&7 &&& &  \\ 2&3 & 4 & 6&7 & 7 &&&&&&& \end{ytableau}, \]
which is partition-shaped, hence highest weight.

\begin{thm}[\LR\ Rule for Ptableaux] \label{LR for ptab}Let ${\mathbb B}^{\mu}(T_{\mmax})$ and ${\mathbb B}^{\nu}(U)$ be two irreducible crystals in $\PTab_{n}$ of highest weights $\mu$ and $\nu$, with highest weight ptableaux $T_{\mmax}$ and $U$, respectively. Then
\[ c_{\mu \nu}^{\lambda} =  ^{\#}\!\!\{ T \in {\mathbb B}^{\nu}(U) :\hbox{$T_{\mmax} \otimes T  $ is partition-shaped, of shape $\lambda$} \}. \]
In particular, the ptableaux $T_{\mmax} \otimes T$ are the highest weight elements in the irreducible constituents of ${\mathbb B}^{\mu}(T_{\mmax})\otimes {\mathbb B}^{\nu}(U)$ of highest weight $\lambda$.\label{LR rule}
\end{thm}

\subsection{Connections with the Classical \LR\ Rule}
Theorem~\ref{LR rule} should be viewed as a generalization of the classical Littlewood-Richardson rule for enumerating \LR\ coefficients $c_{\mu \nu}^{\lambda}$. Below, we make this connection precise, and in doing so, highlight an advantage that the ptableaux model has over the SSYT model for crystals in $[n]^{\otimes k}$.

Let $\mathbb B_{\mu}$ and $\mathbb B_{\nu}$ denote two irreducible crystals of SSYT, of highest weight $\mu$ and $\nu$, respectively. In the literature one finds expressions such as:
$\mathbb B_{\mu} \otimes \mathbb B_{\nu} \cong \bigoplus_{\lambda} c_{\mu \nu}^{\lambda} \mathbb{B_{\lambda}} $
which denote decompositions of a tensor product into irreducible constituents (up to isomorphism). However, as noted above, there is \emph{no} actual definition of a tensor product of SSYT crystals. At best, we \emph{identify} each factor of $\mathbb B_{\mu} \otimes \mathbb B_{\nu}$ with their row reading image (a canonical crystal of words) and compute the tensor product by concatenation of words. But, by virtue of this, we must replace the actual decomposition of $\mathbb B_{\nu} \otimes \mathbb B_{\mu}$ into irreducible constituents with isomorphic copies in the sum $\oplus_{\lambda} c_{\mu \nu}^{\lambda} \mathbb B_{\lambda}$ (the use of ``$\cong$" is essential). That is, typically \emph{no} crystal in $\mathbb B_{\mu} \otimes \mathbb B_{\nu}$ equals $\mathbb B_{\lambda}$, even when identifying it with a crystal of words.

This is not problematic in itself and, indeed, the point of the \LR\ rule is to count multiplicities of isomorphic crystals. The trouble with the SSYT model is that it is of no help in evaluating $c_{\mu \nu}^{\lambda}$, or analyzing its properties. Thus, one needs an entirely \emph{new} combinatorial model, namely \LR\ fillings.

The \LR\ Rule enumerates $c_{\mu \nu}^{\lambda}$ by counting the set of \emph{\LR\ fillings} of the skew shape $\lambda / \mu$ with content $\nu$. To find such a filling, remove the partition $\mu$ from the partition $\lambda$ leaving the skew shape $\lambda / \mu$. In the remaining boxes, using the partition $\nu=(\nu_{1}, \ldots , \nu_{n})$, we place in the skew shape $\nu_{1}$-many $1$'s, $\nu_{2}$-many $2$'s, etc., filling the entire skew shape, so that the filling is semistandard, and such that the filling satisfies the \emph{word condition}, meaning that the row reading word of the filling (in the skew shape) is a Yamanouchi word. The number of such fillings equals $c_{\mu \nu}^{\lambda}$.

We will demonstrate how this classical construction is a special case of our Theorem~\ref{LR rule} above. Recall the \emph{row reading} we adopted for SSYT. We obtain a word $\omega$ from a SSYT $P$ be writing the content in a row from right to left, from the top row to the bottom.

\begin{df} We say a ptableau $T \in \PTab_{n}$ satisfies the \emph{word condition on ptableaux} if and only if $T=\Pf(\omega)$, where $\omega$ is the row-reading word of a SSYT $P$. \label{word df}
\end{df}

We use the term ``word condition" out of tradition as it is used in the context of the Littlewood-Richardson rule. However, the word condition on ptableaux does \emph{not} correspond to the Yamanouchi condition on words in $[n]^{\otimes k}$. That is, if a ptableau $T$ satisfies the word condition, this does \emph{not} imply that $T = \Pf (\omega)$ for some Yamanouchi word $\omega$. Further, while Yamanouchi words (being highest weight elements in $[n]^{\otimes k}$) map to highest weight ptableaux (partition-shaped by Theorem~\ref{highest weight}), when we say a ptableaux \emph{itself} satisfies the word condition, it \emph{need not} be highest weight, and need only satisfy the classical ``word condition'' on \LR\ fillings.

\begin{lem} Suppose $T$ satisfies the word condition for ptableaux. Then $T$ lies in a connected crystal graph such that all ptableaux in the graph satisfy the word condition.
\end{lem}\label{wd cond}
\begin{proof} This follows by noting that the maps fro SSYT to words, and then from words to ptableaux, are isomorphisms of crystal graphs~\cite{BumpSchilling}.
\end{proof}

The following is easily shown:

\begin{lem} A ptableau $T$ satisfies the word condition on ptableaux (Definition~\ref{wd cond}) if and only if the number of $i$'s appearing in $T$ in rows $i$ through $i+k$ is greater than or equal to the number of $(i+1)$'s appearing in rows $(i+1)$ through $i+k+1$, for $k \geq 0$.
\end{lem}\label{wd cond converse}

 Let $\nu = ( \nu_{1}, \ldots , \nu_{k})$ be a partition, and let $T_{\nu}$ denote the ptableau whose first row contains $\nu_{1}$-many $1$'s, whose second row contains $\nu_{2}$-many $2$, etc. So, for example, if $\nu = (4,3,3,1)$ then:
 \ytableausetup{smalltableaux}
\[ T_{\nu} = \begin{ytableau} 1 & 1 & 1 & 1 \\ 2 & 2 &2 & \\ 3 & 3 & 3 & \\ 4 & & & \end{ytableau} . \]

\begin{cor}
If a ptableau $T \in \PTab_{n}$ satisfies the word condition, it lies in some ${\mathbb B}^{\nu}(T_{\nu})$, the crystal with highest weight $T_{\nu}$, where $\nu = (\nu_1, \nu_2, \ldots)$, and $\nu_{1}$ is the number of $1$'s appearing in $T$, $\nu_{2}$ is the number of $2$'s appearing in $T$, etc. In particular, the irreducible crystal ${\mathbb B}^{\nu}(T_{\nu})$ is composed of ptableaux satisfying the word condition on ptableaux. \label{word ptab}
\end{cor}

\begin{proof}
Suppose $T$ is right-justified (the word condition is clearly invariant under row equivalence). The word condition implies that we \emph{only} find $1$'s in row 1 of $T$, only $1$'s and $2$'s in rows 1 and 2 of $T$, etc. Thus, in row $2$, by the word condition the $2$'s in row 2 do not extend left past the $1$'s in row 1, and since $T$ is right-justified, the $1$'s in rows 1 and 2 have no gaps in the columns in which they appear. But continuing, since the $3$'s in row 3 cannot extend to the left past the $2$'s in row 2, then the $2$'s in row 3 cannot extend to the left past the $1$'s in row 2, so that the $1$'s appearing in rows 1,2, and 3 appear in consecutive columns, etc. Thus, the horizontal strip of $1$'s is maximal, in that there is no column in $T$ not containing a $1$. But then those $1$'s may only move directly upward (under the action of some crystal operator $e_{i}$), so that in the highest weight ptableau of $T$, the top row consists entirely of $1$'s. After performing these operations (moving $1$'s upward), we argue inductively to conclude the highest weight is of the form $T_{\nu}$.

Since $T_{\nu}$ satisfies the word condition, then by Lemma~\ref{wd cond} ${\mathbb B}^{\nu}(T_{\nu})$ is composed of ptableaux satisfying the word condition on ptableaux.
\end{proof}

With these results, we see that the classical \LR\ rule is a special case of Theorem~\ref{LR rule} above. Consider some \LR\ filling of skew shape $\lambda / \mu$:

\[ \begin{ytableau} \ & &&&& &1&1&1&1 \\
 & & &1 & 1 & 2 & 2& 2 \\
&1 & 2 & 3 & 3 \\
2 & 3 &  4 & 4 \end{ytableau}.  \]

Here, $\mu = (6,3,1,0)$, $\lambda = (10, 8, 5, 4)$, and $\nu = (7, 5,3,2)$. We can fill the blanks of shape $\mu$ with any semistandard content; for now we use $T_{\mu}$, and distinguish it by circling the content:

\[ \begin{ytableau} \textcircled{$\scriptstyle 1$} &\textcircled{$\scriptstyle 1$} &\textcircled{$\scriptstyle 1$}&\textcircled{$\scriptstyle 1$}&\textcircled{$\scriptstyle 1$}& \textcircled{$\scriptstyle 1$}&1&1&1&1 \\
\textcircled{$\scriptstyle 2$} &\textcircled{$\scriptstyle 2$} &\textcircled{$\scriptstyle 2$} &1 & 1 & 2 & 2& 2 \\
\textcircled{$\scriptstyle 3$}&1 & 2 & 3 & 3 \\
2 & 3 &  4 & 4 \end{ytableau}.  \]

Viewing the above as a (highest weight) \emph{perforated tableau}, we note that, up to a relabeling of content (we would have to add $3$ to each entry of $T$ under our definition of tensor product), we may regard this ptableau as a tensor product:

 \[ \begin{ytableau} \textcircled{$\scriptstyle 1$} &\textcircled{$\scriptstyle 1$} &\textcircled{$\scriptstyle 1$}&\textcircled{$\scriptstyle 1$}&\textcircled{$\scriptstyle 1$}& \textcircled{$\scriptstyle 1$}&1&1&1&1 \\
\textcircled{$\scriptstyle 2$} &\textcircled{$\scriptstyle 2$} &\textcircled{$\scriptstyle 2$} &1 & 1 & 2 & 2& 2 & &\\
\textcircled{$\scriptstyle 3$}&1 & 2 & 3 & 3  & &&&&\\
2 & 3 &  4 & 4&&&&&& \end{ytableau} \rightsquigarrow \begin{ytableau} \textcircled{$\scriptstyle 1$} &\textcircled{$\scriptstyle 1$} &\textcircled{$\scriptstyle 1$}&\textcircled{$\scriptstyle 1$}&\textcircled{$\scriptstyle 1$}& \textcircled{$\scriptstyle 1$}&4&4&4&4 \\
\textcircled{$\scriptstyle 2$} &\textcircled{$\scriptstyle 2$} &\textcircled{$\scriptstyle 2$} &4 & 4 & 5 & 5& 5 & &\\
\textcircled{$\scriptstyle 3$}&4 & 5 & 6 & 6  & &&&&\\
5 & 6 &  7 & 7&&&&&& \end{ytableau} =
\begin{ytableau} \textcircled{$\scriptstyle 1$} &\textcircled{$\scriptstyle 1$} &\textcircled{$\scriptstyle 1$}&\textcircled{$\scriptstyle 1$}&\textcircled{$\scriptstyle 1$}& \textcircled{$\scriptstyle 1$}\\
\textcircled{$\scriptstyle 2$} &\textcircled{$\scriptstyle 2$} &\textcircled{$\scriptstyle 2$}& && \\
\textcircled{$\scriptstyle 3$}&&&&&\\ &&&&&
 \end{ytableau}
\otimes
\begin{ytableau} \ & & &1&1&1&1 \\
 &1 & 1 & 2 & 2& 2 & \\
1 & 2 & 3 & 3 & & & \\
2 & 3 &  4 & 4 & & & \end{ytableau} .  \]

That is, the filling of the skew shape $\lambda / \mu$ is determined by the ptableau
 \[ T = \begin{ytableau} \ & & &1&1&1&1 \\
 &1 & 1 & 2 & 2& 2 & \\
1 & 2 & 3 & 3 & & & \\
2 & 3 &  4 & 4 & & & \end{ytableau}. \]

Since $T$ satisfies the word condition by Lemma~\ref{wd cond converse}, then, by Corollary~\ref{word ptab}, we see that $T$ is in the crystal whose highest weight is $T_{\nu}$.  Thus, by Theorem~\ref{LR rule} the number of $T \in {\mathbb B}^{\nu}(T_{\nu})$ such that $T_{\mu} \otimes T$ is partition shaped (with left-justified content of shape $\lambda$) equals the number of \LR\ fillings of $\lambda / \mu$ of content $\nu$, and each filling is, up to a relabeling, row equivalent to some such ptableau $T$.

\section{Evacuation, Involutions and Highest Weights in $\PTab$}
We now present combinatorial algorithms for computing with highest weight elements in ptableaux crystals. These results show how to connect crystal operators on perforated tableaux to well-known algorithms, including Sch\"{u}tzenberger evacuation and the tableau switching of Benkart, Sottile, and Stroomer (see~\cite{Sottile}). The basis of these connections will be Lemma~\ref{evac ptab} below, showing that the Sch\"{u}tzenberger evacuation map can be ``factored" as a product of ptableaux crystal operators.

\subsection{Evacuation and SSYT}\label{evac intro}

Crystal graphs enjoy a number of symmetries. One that has been extensively studied is the \emph{Lusztig Involution}, denoted here $\phi_{Lus}$, defined by the condition that if $c_{max}$ is the highest weight element in some irreducible crystal graph $C$, and $c = f_{i_{1}} \cdots f_{i_{k}} c_{max}$, then the image of the Lusztig Involution is
\[ \phi_{\Lus}(c)= e_{n-i_{1}} \cdots e_{n-i_{k}} c_{min}, \]
where $c_{min} \in C$ denotes the corresponding lowest element of the crystal.

The Lusztig involution is \emph{defined} by the relations above. However, when nodes of crystal graphs are determined by combinatorial models, various combinatorial algorithms for \emph{computing} $\phi_{Lus}$ have been of interest.

In the SSYT model of crystals, the Lustig image of a SSYT $T$, $\phi_{\Lus}(T)$, is computed by the \emph{Sch\"utzenberger involution} on SSYT (see~\cite{BZ}). The Sch\"utzenberger involution has two steps. The first is the \emph{evacuation} procedure (described below) producing an anti-partition shaped tableau. This is followed by \emph{rotation}, which is the $\Rot$ map defined above, applied to the output of the evacuation step, regarding the output as a perforated tableau (but whose left-justified image under $\Rot$ becomes a SSYT once again).

The evacuation step of the Sch\"utzenberger involution is performed with a series of \emph{jeu de taquin} slides which pass a blank through a given SSYT $T$. Since performing jeu de taquin slides on a ptableau produces another ptableau, while the image of a semistandard tableau after jeu de taquin need not be semistandard (in that it may have interior blanks), it will be convenient to regard $T$ (temporarily) as a left-justified ptableau (so $T = \,^*T$) during the evacuation step.

\begin{df} (See~\cite{lenartI,lenartII}) Define $\phi_{\Lus}$ on crystals of SSYT as follows: First, define the map $\Evac: \PTab_{n} \rightarrow \PTab_{n}$. Given a SSYT $T$, regarded as a left-justified ptableau, an \emph{outer corner} of $T$ is a blank appearing in $^*T$ with no other blanks immediately to its left, or immediately above it, and all blanks below it and to its right. Thus, in $^*T$, an outer corner will have content above it (labeled $b$ below), and to its left (labeled $c$ below):
\[ \begin{ytableau}
\scriptstyle{\dots} &\star & a & b & \scriptstyle{\star}&\scriptstyle{\dots}  \\
\scriptstyle{\dots} &\star  & c & & \scriptstyle{\dots} &  \scriptstyle{\dots} \end{ytableau}. \]

Fix some choice of outer corner, and call this the \emph{distinguished blank}. Then perform \emph{inward jeu de taquin} on this blank as follows: Let $c$ be the content to the left of the distinguished blank (as in the figure above), and $b$ be the content immediately above it.  Iff $b \geq c$, swap the distinguished blank with $b$. If not, (so $b < c$), swap the blank with $c$. Either of these swaps is called a \emph{jeu de taquin slide}, or just a \emph{slide}.
At this point, the distinguished blank will either:

\begin{enumerate}
\item Have content above it and to its left, in which case we again perform a jeu de taquin slide on the distinguished blank.
\item Have content \emph{only} above it, or \emph{only} to its left, in which case we swap the blank with this single content entry.
\item Have only blanks to its left and above it, in which case we do nothing, and consider the blank ``fixed".
\end{enumerate}

We repeat this process, performing slides on the distinguished blank until it has been fixed. We then identify a new outer corner as the distinguished blank, and perform inward jeu de taquin on it. This is repeated until there are no outer corners remaining, resulting in a ptableau $\Evac(T)$. To complete the Sch\"utzenberger involution, we then apply the map $\Rot$ to $\Evac(T)$:
\[ \phi_{\Lus}(T) = \Rot(\Evac(T)). \]
 \end{df}

Berenstein and Zelevinsky (see~\cite{BZ}) proved $\phi_{\Lus}(T) = \Rot \circ \Evac(T)$ on the irreducible crystal containing the SSYT $T$. They also showed that $\Evac(T)$ is independent of the sequence of outer corners chosen in $T$.

The evacuation portion of the algorithm above starts with an anti-partition shape of blanks in the lower-right corner of the diagram that migrate upwards, under inward jeu de taquin, to form a partition-shaped collection of blanks in the upper left, of the same corresponding shape. Equivalently, if a (partition-shaped) ptableau $T$ has weight $(a_{1}, \ldots, a_{k})$ (where $a_{i} \geq a_{i+1}$ since it is partition shaped), then the ptableau $\Evac(T)$ will have weight $(a_{k}, a_{k-1}, \ldots , a_{1})$.

 For example, given a semistandard Young tableau $T$:
\ytableausetup{smalltableaux}
 \[  \begin{ytableau}
1 & 1 & 2 & 2 & 3 &4 \\
2 & 3 & 3 & 4& & \\
3 &4 & 5 & & & \end{ytableau}, \]
we compute $\Evac(T)$ by first performing inward jeu de taquin on the outer corner in row 2:
 \[  \begin{ytableau}
1 & 1 & 2 & 2 & 3 &4 \\
2 & 3 & 3 & 4& & \\
3 &4 & 5 & & & \end{ytableau} \Rightarrow \begin{ytableau}
1 & 1 & 2 & 2 & 3 &4 \\
2 & 3 & 3 & &4 & \\
3 &4 & 5 & & & \end{ytableau} \Rightarrow \begin{ytableau}
1 & 1 & 2 & 2 & 3 &4 \\
2 & 3 &   & 3  &4  & \\
3 &4 & 5 &     &    & \end{ytableau} \Rightarrow \begin{ytableau}
1 & 1 & 2 & 2 & 3 &4 \\
2 &   & 3  & 3  &4  & \\
3 &4 & 5 &     &    & \end{ytableau} \Rightarrow \begin{ytableau}
1 & 1 & 2 & 2 & 3 &4 \\
 &  2 & 3  & 3  &4  & \\
3 &4 & 5 &     &    & \end{ytableau}
 \Rightarrow \begin{ytableau}
 \ & 1 & 2 & 2 & 3 &4 \\
 1& 2  &3  & 3  &4  & \\
3 &4 & 5 &     &    & \end{ytableau}. \]
Continuing inner jeu de taquin, we obtain:

\[ \Evac(T) =\begin{ytableau}
 \ &    &  &  1 &  2 & 3 \\
   &   &   2  &  2  &  3  &4 \\
1 &  3& 3 &  4  & 4 & 5 \end{ytableau}. \]
Thus, (see \cite{BZ}) the composition of $\Evac$ with $\Rot$:

\[ \Rot(\Evac(T)) = \Rot \left(\,  \begin{ytableau}
 \ &    &  & 2  &  2 & 3 \\
   &   &  1   & 2   &  3  &4 \\
1 &  3& 3 &  4  & 4 & 5 \end{ytableau}\, \right) =  \begin{ytableau}
 1 &  2  &  2& 3  &  3 & 5 \\
  2 &  3 &  4   & 5  &    & \\
3 &  4& 4 &     &   &   \end{ytableau}, \]
 results in the SSYT $\phi_{Lus}(T)$, the value of the Lusztig involution in the crystal of SSYT applied to $T$.

 Recall that an arbitrary SSYT $T$ (which can represent an arbitrary node in a crystal of SSYT) may also represent the (left-justified) highest weight node $T$ in some crystal of ptableaux. In crystals of SSYT we can factor the Lusztig involution as $\phi_{Lus} = \Rot \circ \Evac$. In the next section we show that the map $\Evac$ itself factors as a product of ptableaux crystal operators. In fact, we show that, regarding some partition-shaped $T$ as a highest weight ptableau, the map $\Evac$ takes $T$ to its associated \emph{lowest weight} ptableau.
 In the example above one can check $f_{1}f_{2}f_{2}f_{2}f_{1}f_{1}T = \Evac(T)$. Since $\Evac(T)$ is anti-partition shaped and is the image of $T$ under crystal operators, it is the lowest weight element in the irreducible crystal with highest weight $T$, by Theorem~\ref{highest weight}. It can be shown that the sequence of indices of the $f_{i}$'s used in this process form a ``good word" in the language of string patterns (see~\cite{BumpSchilling}).

Thus, the evacuation map gives a combinatorial algorithm to compute lowest weight ptableaux from their associated highest weights, as well as interpreting jeu de taquin slides as crystal operators on perforated tableaux.
 Finding a combinatorial algorithm to compute Lusztig involutions for \emph{arbitrary} ptableaux is outside the scope of this paper, but can be achieved (and requires an analysis of the Robinson-Schensted construction in the setting of ptableaux~\cite{A-W RSK}).

\subsection{Evacuation and Crystal Operators on Ptableaux}

\begin{df} We say $T \in \PTab$ is \emph{skew-shaped} if it is row-equivalent to a ptableau in which the only blanks appearing are a partition-shaped collection of blanks in the upper left corner, and an anti-partition-shaped collection of blanks in its lower right corner. We call this the \emph{skew-shaped form} of $T$.

We say a blank of a skew-shaped ptableau $T$ is an \emph{inner corner} of the skew-shaped form of $T$ if it appears in the partition-shaped collection of blanks in the upper left-corner, and if there are no blanks to its right or below it.
\end{df}
In other words, a ptableau is skew-shaped if, in its skew-shaped form, the content is skew-shaped in the traditional sense of tableaux.

Above we made use of \emph{inward jeu de taquin} in the definition of the map $\Evac$, where a blank in an outer corner is moved left or upward by jeu de taquin slides. Below we consider \emph{outward jeu de taquin}, where blanks starting in an inner corner are swapped with content below or to the right of the blank.

\begin{lem} \label{e evac} Suppose $T$ is a skew-shaped ptableau. Choose an \emph{inner corner} $S$ of this skew shape (that is, $S$ is a blank in the skew-shaped form of $T$ that has content directly to its right and directly below it), appearing in row $s$. We perform outward jeu de taquin on this skew-shaped form, starting at $S$, and call the resulting skew-shaped ptableau $T'$. Let ${\cal P}$ denote the outward jeu de taquin path of $S$ from the inner corner of $T$ in row $s$ to an outer corner of $T'$ appearing in some row $t+1$, where $s \leq t$. Then, up to row-equivalence, we have
\[ T' = e_{t} e_{t-1} \cdots e_{s} T, \] \label{evac ptab}
and the path of $S$, resulting from these crystal operators, is also ${\cal P}$.
\end{lem}
\begin{proof}
Let $T$ be a skew-shaped ptableau, and for simplicty, denote its (fixed) skew-shape form by $T$. Fix an inner corner of $T$, denoted $\begin{ytableau}\circ\end{ytableau}$, appearing in row $s$ (shaded gray below). To perform outward jeu de taquin, we swap
$\begin{ytableau}\circ\end{ytableau}$ with content in $T$ as long as all entries remain column strict.

We prove the claim inductively by showing that, for each step of the outward jeu de taquin path of $\begin{ytableau}\circ\end{ytableau}$, if the slide moves $\begin{ytableau}\circ\end{ytableau}$ from row $k$ to row $k+1$, the same (row-equivalent) ptableau is obtained by applying $e_k$ to the ptableau determined prior to that slide.

Suppose we performed a sequence of jeu de taquin slides, stopping when $\begin{ytableau}\circ\end{ytableau}$ appears in row $i$, and is no longer able to move right without creating a column strictness violation. Call the resulting ptableau $T'$ (where the $\begin{ytableau}\star\end{ytableau}$ denotes content):


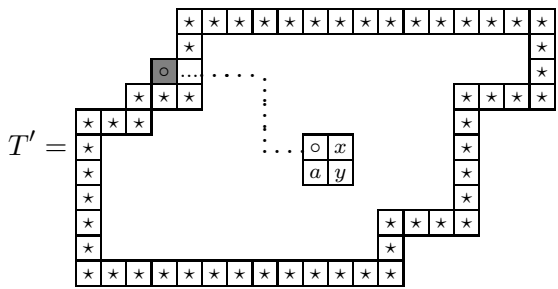
\[
T' =\begin{ytableau} \none & \none& \none & \none &\star &\star &\star &\star & \star&\star &\star & \star &\star &\star &  \star &\star &\star&\star&\star\\
\none & \none& \none & \none & \star&\none &\none &\none &\none &\none &\none &\none &\none &\none &\none &\none&\none&\none& \star \\
\none & \none& \none & *(gray)\circ & \dots&\none[\dots] &\none[\dots] &\none &\none &\none &\none &\none &\none &\none &\none &\none & \none & \none&  \star\\
\none& \none& \star  &  \star & \star &\none &\none &\none[\vdots] &\none &\none &\none &\none &\none &\none &\none &  \star& \star& \star& \star\\
 \star & \star & \star  & \none &\none &\none &\none &\none[\vdots] &\none &\none &\none &\none &\none &\none &\none & \star \\
 \star & \none& \none & \none &\none &\none &\none &\none[\vdots] &\none[\dots]  &\circ &x &\none &\none &\none &\none &  \star\\
   \star& \none& \none & \none &\none &\none &\none &\none &\none  &a &y &\none &\none &\none &\none & \star \\
    \star& \none& \none & \none &\none &\none &\none &\none &\none  &\none &\none &\none &\none &\none &\none & \star \\
  \star & \none& \none & \none &\none &\none &\none &\none &\none  &\none &\none &\none &  \star& \star & \star& \star \\
 \star  & \none& \none & \none &\none &\none &\none &\none &\none  &\none &\none &\none &  \star&\none &\none&\none \\
  \star &  \star& \star  & \star  &  \star& \star & \star & \star &  \star&  \star& \star & \star & \star &\none &\none&\none \\
\end{ytableau} \]

We assume, inductively, that this ptableau is row-equivalent to $e_{i-1} \cdots e_s T$.
Because the blank \framebox{$\circ$} cannot move farther right at this point, we must be unable to move the box \framebox{$x$} over the box \framebox{$a$}, and so $a \leq x$.

The evacuation procedure thus swaps the blank \framebox{$\circ$} in row $i$ with the content \framebox{$a$} in row $(i+1)$.
We show that, after we left justify $T'$, computing $e_{i}(T')$ will swap \framebox{$a$} with the blank directly above it, which will show that the evacuation procedure matches the crystal operators as described in the statement of the lemma.

Since $T$ is initially skew-shaped, we may assume $T[i,i+1]$, before we left justified, had the form:
\[ \begin{ytableau} \none & \none & \cdots  &\cdots & \circ & x & \cdots & \cdots & \cdots & \cdots \\ \cdots & \cdots &\cdots & \cdots & a & y & \cdots \end{ytableau}, \]
where the boxes to the left of $\framebox x$ are filled.

Thus, in order to complete the left-justification of $T[i,i+1]$ we need only move the content to the \emph{left} of the blank $\framebox{$\circ$}$ in row $i$  (possibly) further left, but all other content in row $i$ cannot move. Thus, the blank $\framebox{$\circ$}$ is necessarily the rightmost blank in row $i$ over content in row $i+1$, so by the definition of $e_i$ we move the entry $\framebox{a}$ upward into row $i$. If evacuation swaps $\framebox{$\circ$}$ down and the content $\framebox{$a$}$ up, then upon left justification, $e_{i}(T)$ would swap $\framebox{$a$}$ up as well.
\end{proof}

By duality under the $\Rot$ map, we then have:

\begin{cor}
\label{f evac} Suppose a ptableau $T$ has a skew shape row-equivalent form. Choose an outer corner $S$ of this skew shape (so $S$ is a blank that has content directly to its left and directly above it), appearing in row $j+1$, and perform inward jeu de taquin, moving $S$ from the lower outer corner of $T$ in row $j+1$ to an inner corner of $T'$ in row $i$, where $i \leq j$. Then, up to row equivalence, we have
\[ T' = f_{i} f_{i-1} \cdots f_{j} T. \]
\end{cor}
We now prove the results mentioned at the end of Section~\ref{evac intro}.

\begin{thm} \begin{enumerate}
\item Let $T \in \SSYT(\nu)$ be a semistandard Young tableau of shape $\nu$, or equivalently, a partition-shaped ptableau of weight $\nu$. Then, there is a sequence of ptableaux crystal operators $ f_{i_{1}}, \ldots  , f_{i_{s}}$ such that
\[ Evac ( T )  =f_{i_{1}}, \ldots  , f_{i_{s}} T. \]
In particular, each evacuation path constructed in computing $Evac(T)$ by inward jeu de taquin may be ``factored'' as a product of ptableaux crystal operators. \label{involution}
\item If $T_{\min}$ denotes the lowest weight element in the ptableaux crystal containing the (partition-shaped) highest weight element $T$, then
\[ T_{\min} = \Evac(T). \]
\end{enumerate}
\end{thm}
 \begin{proof} We identify some SSYT $T$ as a highest weight ptableau $T=T_{max} \in \PTab_n$, by Theorem~\ref{highest weight}. We assume, for now, that $n$ is chosen so that there are no entirely blank rows in $T_{\mmax}$. Thus, in its left-justified form, all blanks in $T_{\mmax}$ appear in the lower right corner, and this collection of blanks is anti-partition shaped (the number of blanks weakly increases as we proceed down the rows).

 Applying $\Evac$ to $T_{\mmax}$ moves these blanks to the upper left corner. We choose to move one blank entirely through $T_{\mmax}$ before staring another, so that Corollary~\ref{f evac} applies. Once all blanks in the anti-partition shape of blanks below $T_{\mmax}$ have been moved above it, $T_{\mmax}$ has been transformed into the form of a lowest weight ptableau (which we denote $T_{\mmin}$) by Theorem~\ref{highest weight}. By Corollary~\ref{f evac}, this ptableau is in the same crystal as $T_{\mmax}$.

 Thus, as in the example above, the transformation from the highest weight ptableau

  \[  T_{\mmax} = \begin{ytableau}
1 & 1 & 2 & 2 & 3 &4 \\
2 & 3 & 3 & 4& & \\
3 &4 & 5 & & & \end{ytableau}, \,\,\hbox{to the lowest weight,} \,\,
 T_{\mmin}= \begin{ytableau}
 \ &    &  &  1 &  2 & 3 \\
   &   &   2  &  2  &  3  &4 \\
1 &  3& 3 &  4  & 4 & 5 \end{ytableau}, \]
by evacuation is also obtained by an appropriate sequence of crystal operators $f_{i}$.

 For the general case, we argue by example (the general principle will be clear). Assume that the highest weight ptableau $T_{\mmax}$ has some rows of blanks:
  \[  T_{\mmax} = \begin{ytableau}
1 & 1 & 2 & 2 & 3 &4 \\
2 & 3 & 3 & 4& & \\
3 &4 & 5 & & &  \\
&&&&& \\
&&&&& \\ &&&&& \end{ytableau}. \]

Applying evacuation in the first three rows of $T_{\mmax}$, or applying the same appropriate set of crystal operators $f_{i}$, will result in the ptableau:

 \[  T' = \begin{ytableau}
\ &    &  &  1 &  2 & 3 \\
   &   &   2  &  2  &  3  &4 \\
1 &  3& 3 &  4  & 4 & 5 \\
&&&&& \\
&&&&& \\ &&&&& \end{ytableau}\, .
\quad \hbox{But then} \quad
( f_{1})^{3}(f_{2})^{4}(f_{3})^{6} (T') = \begin{ytableau}
\ &&&&& \\
 &    &  &  1 &  2 & 3 \\
   &   &   2  &  2  &  3  &4 \\
1 &  3& 3 &  4  & 4 & 5 \\
&&&&& \\ &&&&& \end{ytableau}\, .\]
We now repeat the process to move blanks in the bottom two rows upward, either by evacuation or crystal operators, finally obtaining the lowest weight:
\[  T_{\mmin} = \begin{ytableau}
\ &&&&& \\
&&&&& \\ &&&&& \\
 &    &  &  1 &  2 & 3 \\
   &   &   2  &  2  &  3  &4 \\
1 &  3& 3 &  4  & 4 & 5
 \end{ytableau}, \]
and the general result follows.

The second claim of the Theorem follows by the second claim of Theorem~\ref{highest weight} since $\Evac(T_{\mmax})$ is anti-partition shaped and in the image of $T_{\mmax}$ under crystal operators. \end{proof}

We should note that, as pointed out in~\cite{Shimozono}, computing commutators for tensor products of words can be reduced, algorithmically, to computing commutators of individual factors, and that, in turn, this computation may be computed by \emph{jeu de taquin} slides on tableaux formed by those factors. Further, subsequent to our results the recent work of Gerber and Lecouvey~\cite{GerberLecouvey} has appeared and shows that \emph{jeu de taquin} slides on tableaux can be related to crystal operators on the corresponding duals of those objects. However, as pointed out earlier, these results do not (quite) imply ours, as we show, in the case of ptableaux, that the \emph{jeu de taquin} moves are \emph{themselves} the images of crystal operators on the \emph{same} object, and not its dual.

\subsection{Commutators and Highest Weights}
We now consider the problem of combinatorially determining highest weight ptableaux in tensor products of crystals, using some of the tools of the previous section.  Let $B_{\mu}$ and $B_{\nu}$ be two irreducible crystals of ptableaux, with highest weight elements of weight $\mu$ and $\nu$ respectively. This tensor product decomposes:
\[ B_{\mu} \otimes B_{\nu} \cong \bigoplus_{\lambda} B_{\lambda}^{\oplus c_{\mu \nu}^{\lambda}}, \]
where, in the the \emph{isomorphism} above, $B_{\lambda}$ is a \emph{fixed} irreducible crystal of highest weight $\lambda$, and the \LR\ coefficient, $c_{\mu \nu}^{\lambda}$ counts its \emph{multiplicity} in the decomposition. However, $B_{\mu}$ and $B_{\nu}$ denote \emph{specific} crystals of ptableaux (corresponding to specific words in $[n]^{\otimes k}$), so, in addition to counting the multiplicity of irreducible crystals in a given isomorphism type, we should try to distinguish and compute \emph{within} a particular, realized tensor product.

By Theorem~\ref{LR for ptab} we know that a highest weight element of an irreducible constituent of the tensor product $B_{\mu} \otimes B_{\nu}$ has the form $T_{\mu\,\mmax} \otimes T$ where $T_{\mu\,\mmax}$ is the highest weight element of $B_{\mu}$, $T \in B_{\nu}$, and $T_{\mu\,\mmax } \otimes T$ is a partition-shaped ptableau, of some shape $\lambda$ (so $c_{\mu \nu}^{\lambda} \neq 0$).

Further, it is well known that the crystals $B_{\mu} \otimes B_{\nu}$ and $B_{\nu} \otimes B_{\mu}$ are isomorphic. A crystal isomorphism $\tau : B_{\mu} \otimes B_{\nu} \rightarrow B_{\nu} \otimes B_{\mu}$ is called a ``commutator". Given some element $T \otimes T' \in B_{\mu} \otimes B_{\nu}$, Henriques and Kamnitzer~\cite{kam} compute a commutator $\tau(T \otimes T') \in B_{\nu} \otimes B_{\mu}$ by means of the Luztig involution, $\tau(T \otimes T') = \phi_{\Lus}(\phi_{\Lus}(T')\otimes \phi_{\Lus}(T_{\inv}))$. One would hope, however, for a purely combinatorial realization of this map in terms of ptableaux. The general case has been obtained, but is outside the scope of this paper~\cite{A-W RSK}. Here, however, we present a combinatorial algorithm to compute commutators of \emph{highest weight} elements in tensor products, which is of interest not only as a combinatorial algorithm, but in its connections to similar algorithms appearing in other contexts.

A special case of the commutator map in the context of \LR\ fillings (but not crystals) has been known for many years. Recall that (left-justified) highest weight ptableaux can be viewed as SSYT. In particular, classical results on \LR\ fillings can be viewed, using ptableaux, as special tensor products resulting in highest weight ptableaux (by Theorem~\ref{LR rule} above). The classical \LR\ rule was obtained as a special case of this result, when \emph{both} $B_{\mu}$ and $B_{\nu}$ are the unique crystals (of their respective weights) that correspond to reading words of SSYT (see Corollary~\ref{word ptab} and the discussion following).

In this special setting, the commutator is defined as a bijection between \LR\ filling of $\lambda / \mu$ with content $\nu$ (where $\lambda$ is the weight of the fixed, irreducible constituent of $B_{\mu} \otimes B_{\nu}$), and a corresponding \LR\ filling of $\lambda / \nu$ of content $\mu$ (corresponding to constituents of $B_{\nu} \otimes B_{\mu}$). James and Kerber~\cite{J-K} give this map by defining a combinatorial algorithm that moves the content of the \LR\ filling of content $\nu$ ``through" the base of $\mu$, so that the tableau becomes one of base $\nu$, with  content $\mu$.

Benkart, Sottile, and Stroomer~\cite{Sottile} show that James' and Kerber's algorithm is a special case of a more general ``tableau switching" algorithm.

We make use of the work of Benkart, Sottile and Stroommer (BSS) in what follows, and also describe a combinatorial algorithm on ptableaux, in the case of highest weights, that will map a given highest weight element of $B_{\mu} \otimes B_{\nu}$ to a canonical corresponding highest weight ptableau in $B_{\nu} \otimes B_{\mu}$. That is, given some $T_{\mu\,\mmax} \otimes T \in B_{\mu} \otimes B_{\nu}$, where $T_{\mu\,\mmax}$ and $T_{\mu\,\mmax} \otimes T$ are both highest weight, we define a combinatorial algorithm $\tau(T_{\mu\,\mmax} \otimes T) = T_{\nu\,\mmax }\otimes T''$ such that
\begin{enumerate}
\item $T'' \in B_{\mu}$ (the \emph{same} crystal,  not merely one isomorphic to $B_{\mu}$).
\item $T_{\nu\,\mmax} \in B_{\nu}$ is the highest weight element in $B_{\nu}$.
\item $T_{\nu\,\mmax }\otimes T''$ is highest weight, and $\wt(T_{\mu\,\mmax } \otimes T) = \wt(T_{\nu\,\mmax}\otimes T'')$.
\end{enumerate}

The algorithm can be described as a type of tableau switching (in the sense of BSS), an evacuation procedure, or a coordinated set of crystal operators. We use all three interpretations to prove the equivalence of these algorithms, and to show they satisfy the requirements above.

We use Lemma~\ref{e evac} and its Corollary~\ref{f evac}. We also need the lemma below, which will help keep track of pairs of evacuation paths.

\begin{lem} \label{evac paths}
Suppose $\Box_{1}$ and $\Box_{2}$ are outer corners of a skew shaped ptableau $T$, with $\Box_{1}$ weakly below and strictly left of $\Box_{2}$.  Then, the evacuation path of $\Box_{1}$ lies weakly below and strictly left of the evacuation path of $\Box_{2}$.
\end{lem}
\begin{proof}
Note that, in the diagrams below, we label $\Box_1$ and $\Box_2$ as \framebox{$1$} and  \framebox{$2$}.
Assume that at some stage of the evacuation of $\Box_1$, we have moved it as far left as it can go in its current row.  Immediately before moving $\Box_1$ up, $T$ contains
\[ \begin{ytableau}  \none &  b & c\\ a &1 & d\end{ytableau} . \]

 Because we cannot move $\Box_1$ left,  we have that $a\le b$.  After swapping $\Box_1$ and \framebox{$b$}, T has
\[ \begin{ytableau}  \none &1 & c\\ a &b & d\end{ytableau}.  \]

Now assume that at some stage in the evacuation of $\Box_2$, it has moved left as far as possible, to the column immediately to the right of the vertical path of $\Box_1$. So $\Box_2$ occupies the box that contained $d$, or occupies a higher box in that column, with an entry $b'<b$ immediately to its left.  (In reality, $\Box_2$ might be several rows below $\Box_1$, but the argument will still hold.)
 \[ \begin{ytableau}  \none &1& c\\ a &b &2\end{ytableau}  \text{ or } \begin{ytableau}  \none &1 & c'\\ \none &b' & 2\\ a &b & d\end{ytableau}. \]

We must show that $\Box_2$ now moves up, so that its evacuation path does not cross to the left of  $\Box_1$'s evacuation path. In the first case, as noted above, $b\leq c$ and so to maintain $T$'s column-strictness, $\Box_2$ must move up.  In the second case, because $\Box_1$ swapped with $b'$, by column-strictness $b'\le c'$.  Thus to maintain column-strictness, $\Box_2$ must swap with $c'$ above it.  Hence, in either case, $\Box_2$ moves vertically, and does not cross the vertical path of $\Box_1$.

$\Box_1$ started weakly below and strictly to the left of $\Box_2$, and so the evacuation path of $\Box_1$ must be weakly below and strictly left of that of $\Box_2$.  \end{proof}

In the work of Benkart, Sottile and Stroomer (BSS)~\cite{Sottile}, what we define below as \emph {BSS perforated tableaux} are called \emph{perforated tableaux}.  Their definition omits Properties 2 through 5 in our Definition \ref {ptab}:

\begin{df} [See~\cite{Sottile}] \label {BSS}Let $T$ be a tableau of numbers and blanks, with numbers from $[n]$. Then $T$ is a \emph{BSS perforated tableau} if, ignoring blanks, the entries of $T$ are weakly increasing in rows and strictly increasing in columns.
\end{df}

The ptableaux defined here (Definition \ref {ptab}) are BSS perforated, but the converse does not hold. For example, the tableau

\[ \begin{ytableau} \none & \ & 1 & 1 & & 2 \\ 2 & & 3 &  & 4 & \\
& 2 & &3 & &\none \\ & 3 & & 4&\none&\none \end{ytableau} \]
is BSS perforated, but not a ptableau, since the $2$'s, $3$'s, and $4$'s do not form horizontal strips.

\begin{df} Suppose we have a tableau with content of two ``types," for example, $1,2,3,\dots$, and
$\bf1,\bf2,\bf3\dots$,
 but no blanks.   If the portion of the tableau formed by each type of number separately form BSS perforated tableaux, then the entire tableau is called a \emph{BSS perforated pair}.\label{BSS}
\end{df}

For example, the diagram below is a BSS perforated pair:
\[ \begin{ytableau}\bf 1&\bf 1 & 1 & 1 &\bf 2& 2 \\
2 &\bf 2& 3 &\bf 2& 4 &\bf 2\\
\bf 2& 2 &\bf 3&3 &\bf 3& \none\\
\bf 3& 3 &\bf 4& 4&\none&\none \end{ytableau} \]

\noindent {\bfseries{\scshape{Algorithm 1 ``Push Down":}}}  Given: $S= T_{\mu\, \mmax } \otimes T$, where $T_{\mu\, \mmax }$ and $T_{\mu\, \mmax } \otimes T$ are highest weight (partition shaped) ptableaux. Call the filling of $T_{\mu\, \mmax }$ (resp. $T$) the \emph{$\mu$ content} (resp. $\nu$ content) of $S$. At the start of the algorithm, the $\nu$ content forms a skew shape.

Let $k$ be the largest value in the $\mu$ content. The rightmost, highest $k$ in the $\mu$ content lies in an inner corner of the skew shape of the $\nu$ content.

Use outer evacuation to push the highest, rightmost $k$ through the $\nu$ content, until it emerges in an outer corner of the (still skew shaped) $\nu$ content. We work through the $k$-strip of the $\mu$ content, moving right to left, top to bottom, at each stage pushing a single $k$ through the $\nu$ content. At the end of this step, the $\mu$ content of $1$'s through $(k-1)$'s lies in a partition shape above the skew shape of the $\nu$ content, and below the $\nu$ content, there is a single horizontal strip of $k$'s, which is part of the $\mu$ content.

We then repeat this pushing down process on the $(k-1)$ strip of the $\mu$ content, starting with the rightmost, highest $k-1$, then the $(k-2)$ strip, etc. At the end of this process, the $\nu$ content forms a partition shape, with the $\mu$ content forming a skew shape below it.

\bigskip

\noindent {\bfseries{\scshape{Algorithm 2 ``Push Up":}}}
In this case, we begin with the $1$-strip of the $\nu$ content. Beginning with the lowest, left-most 1 in the strip, use upper evacuation to push the $1$ up through the $\mu$ content, until it occupies the top left corner of the diagram. Repeat with the rest of the $1$-strip, working left to right, bottom to top. Do the same for the $2$ strip, etc., until the $\nu$ content forms a partition shape, with the $\mu$ content forming a skew shape below it.

\bigskip

We show below that both the Push Down and the Push Up algorithms result in the \emph{same} resulting output, namely, a ptableau in the tensor product $B_{\nu} \otimes B_{\mu}$ of highest weight, in other words, the commutator of $S$.

\begin{thm} Let $B_{\mu}$ and $B_{\nu}$ denote arbitrary irreducible crystals of ptableaux of highest weights $\mu$ and $\nu$, respectively. Let $T_{\mu \,\mmax} \in B_{\mu}$ be the  highest weight ptableau in $B_{\mu}$, and let $T \in B_{\nu}$ be such that $T_{\mu \,\mmax}\otimes T \in B_{\mu} \otimes B_{\nu}$ is also highest weight. Then, applying either the Push Up or the Push Down algorithm to $T_{\mu \,\mmax}\otimes T $ produces the unique highest weight ptableau $T_{\nu \,\mmax}\otimes T'   \in B_{\nu} \otimes B_{\mu}$.
\end{thm}

\begin{proof}
The Push Down algorithm pushes $\mu$ content through the $\nu$ content by outer evacuation, keeping the $\mu$ content BSS perforated. Lemma~\ref{e evac} shows that pushing $\mu$ content down can be accomplished by performing $e_{i}$ crystal operators on $T$. Thus, not only is the image of $T$ a perforated tableau during the Push Down algorithm, its image lies in the same crystal $B_{\nu}$, and not merely one isomorphic to it. Similarly, the Push Up algorithm moves the $\nu$ content in $T$ so that it is BSS perforated, and will have the effect as applying $f_{i}$ crystal operators on $T_{\mu\,\mmax}$, taking it to some ptableau in the same crystal $B_{\mu}$.

Consider the Push Down algorithm applied to $T_{\mu \,\mmax} \otimes T$. By Lemma~\ref{evac paths}, each box in the $k$-strip of the $\mu$ content will arrive, via outer evacuation, to the right of the boxes of the $k$-strip that came before it. Thus, as we move the $k$-strip down through the $\nu$ content, it forms a horizontal strip and, in particular, the $\mu$ content formed by the Push Down algorithm forms a ptableau. From this, it is clear that at each stage of the Push Down algorithm (including when a box of $\mu$ content is in the interior of the $\nu$ content), the resulting tableau forms a perforated BSS pair. So, the Push Down algorithm (defined via outer evacuation) is equivalent to performing BSS tableau switching.

An analogous argument implies that the Push Up algorithm is also a special case of BSS tableau switching.
However, Theorem 2.2  of ~\cite{Sottile} demonstrates that the resulting tableau obtained by BSS tableau switching is independent of the sequence of switches used. Thus, the outputs of the Push Up and the Push Down algorithms are the same.

Hence, while we start with the tensor product $T_{\mu \,\mmax} \otimes T$, with $T_{\mu \,\mmax}$ highest weight in $B_{\mu}$ and $T \in B_{\nu}$, the resulting ptableau (from either algorithm) is a tensor product $T_{\nu \,\mmax}\otimes T'$ with $T_{\nu \,\mmax}$ highest weight in $B_{\nu}$ and $T' \in B_{\mu}$. Since the action of the Push Algorithms applies crystal operators to either the $\mu$ or the $\nu$ content, and the action of both is a BSS switching, the result of either algorithm is the same.
BSS switching preserves shape, so the resulting ptableau is highest weight of the same weight $\lambda = \wt(T_{\mu \,\mmax}\otimes T)$.
\end{proof}

\section{Questions for Future Work}
\begin{enumerate}
\item{\bf Insertion Algorithms.} There is an elegant version of ``RSK"-style results, but couched in the language of ptableaux. In~\cite{A-W RSK} we show that an insertion algorithm on a ptableau $T$ can produce a pair of ptableaux $(P(T), T_{max})$ such that (1) Given any two ptableaux $T$ and $T'$, if $P(T) = P(T')$, then $T$ and $T'$ are plactically equivalent. (2) $T_{max}$ is the highest weight element in the irreducible crystal graph containing $T$, and (3) There are simple algorithms to convert $P(T)$ to a SSYT $P$ and $T_{max}$ to a SSYT $Q$ such that $(P,Q)$ is the RSK pair associated to the word $\omega$, with $\Pf(\omega) = T$. This is natural in the sense that, for example, instead of using $Q$ as a \emph{proxy} for determining if two ptableaux lie in the same irreducible crystal graph, an insertion algorithm computes their respective highest weights directly.
\item{\bf Commutators.} The results obtained in~\cite{A-W RSK} have been used to give a combinatorial algorithm to compute commutators for arbitrary ptableaux, extending the results on highest weights obtained here.
\item {\bf Other Types.} The successes in type $A_{n-1}$ suggest looking for appropriate perforated tableaux models in other types, by looking at tensor products of the standard crystals there. We have determined ptableaux models in all classical types and are investigating their usefulness in those settings.
\end{enumerate}


\begin{thebibliography}{99}

\bibitem{A-W RSK}G.\ Appleby and T.\ Whitehead. RSK and PTableaux (in preparation), (2022)

\bibitem{A-W Ptab} G.\ Appleby and T.\ Whitehead. Perforated Tableaux: A Combinatorial Model for Crystal Graphs in Type $A_{n-1}$. arXiv:2007.11721 (2021)


\bibitem{Sottile} G.\ Benkart, F.\ Sottile and J.\ Stroomer: Tableau
Switching: Algorithms and Applications. J. Combin. Th. Ser. A, {\bf 76}, 11-43, (1996)

\bibitem{BZ} A. Berenstein and A. Zelevinsky: Canonical bases for the quantum group of type $A_{r}$ and piecewise-linear combinatorics.
Duke Math. J. {\bf 82}, 3, 473-502 (1996)

\bibitem{BumpSchilling} D.\ Bump and A.\ Schilling: Crystal Bases, Representations and Combinatorics. World Scientific, Hackensack, NJ (2017)

\bibitem{FuLascoux} A.\ Fu and A.\ Lascoux: Non-symmetric Cauchy kernels for the classical groups. \emph{JCTA}, 116:903-917 (2009)

\bibitem{Fulton} W.\ Fulton. Young Tableaux. Cambridge University Press (1997)

\bibitem{GerberLecouvey} T.\ Gerber and C.\ Lecouvey. Duality and bicrystals on infinite binary matrices. ArXiv: 2009.10397 (2020)

\bibitem{kam} A.\ Henriques and J.\ Kamnitzer: The Octahedron Recurrence and ${\mathfrak gl}_{n}$ Crystals. Advances in Mathematics, {\bf 206}, 1, 211-249 (2006)

\bibitem{HeoKwon} T.\ Heo and J. Kwon. Combinatorial Howe duality of symplectic type. ArXiv:2008.05093 (2020)


\bibitem{J-K} G. James, A Kerber: Representation Theory of the Symmetric Group. Addison-Wesley, New York (1982)

\bibitem{Kash-Nak} M.\ Kashiwara and T.\ Nakashima: Crystal Graphs for representations of the $q$-analogue of classical Lie algebras, \emph{J. Algebra}, {\bf 165}, 2, pp. 295-345 (1994)

\bibitem{Lee} J.\ Lee. Crystal structure on King tableaux
and semistandard oscillating tableaux. ArXiv:1910.04459 (2019)

\bibitem{lenartI} C. \ Lenart: On the combinatorics of crystal graphs, I. Lusztig's Involution. Advances in Mathematics. {\bf 211}, 1, 204-243 (2007)

\bibitem{lenartII} C. \ Lenart: On the combinatorics of crystal graphs, II. The crystal commutator. Proc. Amer. Math. Soc. {\bf 136} 825-837 (2008)

\bibitem{NakayashikiYamada} A.\ Nakayashiki and Y.\ Yamada. Kostka polynomials and energy functions on solvable lattice models. \emph{Selecta Math.}, 3:547-599 (1997)

\bibitem{Shimozono} M.\ Shimozono. Crystals for Dummies. https://www.aimath.org/WWN/kostka/crysdumb.pdf (2005)

\bibitem{vanLeeuwen} M.\ van Leeuwen. Double Crystals of Binary and Integral Matrices. \emph{Electronic Jounal of Combinatorics}. 13 (2006)


\end{thebibliography}
\end{document}